\documentclass[11pt, leqno]{amsart}
\usepackage{amscd, amsmath, amssymb, amsfonts}
\usepackage{latexsym}
\usepackage{mathrsfs}

\usepackage{fullpage}
\usepackage{microtype}
\usepackage{colonequals}

\usepackage{graphicx}
\usepackage{tikz}
\usepackage{tikz-cd}
\usetikzlibrary{calc, positioning}
\usepackage{circuitikz}
\usepackage{mathbbol}

\usepackage{enumitem} 
\setlist[enumerate]{label=\textnormal{(\arabic*)}}

\ctikzset{bipoles/resistor/height=0.15}
\ctikzset{bipoles/resistor/width=0.4}

\definecolor{cadmiumgreen}{rgb}{0.0, 0.42, 0.24}
\definecolor{darkred}{rgb}{.85,0,0}
\definecolor{aqua}{rgb}{0.0, 1.0, 1.0}
\definecolor{byzant}{rgb}{0.74, 0.2, 0.64}

\usepackage[
colorlinks, 
linkcolor=blue,
urlcolor=darkred,
citecolor=cadmiumgreen,
pdfauthor={Omid Amini, Shu Kawaguchi, JuAe Song}, 
pdftitle={},
pdfstartview ={FitV},
]{hyperref}

\newtheorem{Theorem}{Theorem}[section]
\newtheorem{Proposition}[Theorem]{Proposition}
\newtheorem{Lemma}[Theorem]{Lemma}
\newtheorem{Corollary}[Theorem]{Corollary}

\numberwithin{equation}{section}

\theoremstyle{definition}
\newtheorem{defii}[Theorem]{Definition}
\newenvironment{Definition}
  {\pushQED{\qed}\defii}
  {\popQED\enddefii}

\newtheorem{exx}[Theorem]{Example}
\newenvironment{Example}
  {\pushQED{\qed}\exx}
  {\popQED\endexx}

\newcommand{\TT}{{\mathbb{T}}}

\newcommand{\NN}{{\mathbb{N}}}
\newcommand{\ZZ}{\ss{\mathbb{Z}}}
\newcommand{\QQ}{\ssub{\mathbb{Q}}!}
\newcommand{\RR}{\ssub{\mathbb{R}}!}

\newcommand{\OO}{{\mathcal{O}}}

\newcommand{\Ucal}{{\mathcal{U}}}

\newcommand{\Lscr}{{\mathscr{L}}}

\newcommand{\Sscr}{{\mathscr{S}}}

\newcommand{\Xscr}{{\mathscr{X}}}

\newcommand{\rest}[2]{\left.{#1}\right\vert_{{#2}}}  

\newcommand{\Spec}{{\operatorname{Spec}}}

\newcommand{\trop}{\operatorname{trop}}
\newcommand{\an}{{\operatorname{an}}}

\newcommand{\dprime}{{\prime\prime}}
\newcommand{\suppaux}[2]{\scalebox{1}[1.1]{$#1\lvert$}#2\scalebox{1}[1.1]{$#1\rvert$}}
  \newcommand{\supp}[1]{\mathpalette\suppaux{#1}}

\usepackage{mathrsfs}

\usepackage{xcolor}  

\usepackage{todonotes} 

\newcommand{\rec}{\mathrm{rec}}

\newcommand{\z}{\textsc{z}}
\newcommand{\w}{\textsc{w}}
\newcommand{\x}{\textsc{x}}

\newcommand{\Rat}{\operatorname{Rat}}
\newcommand{\zero}{\operatorname{div}}

\newcommand{\Sk}{\operatorname{Sk}}

\newcommand{\f}{{\textnormal{\textsl{\textsf f}}}} 
\newcommand{\e}{{\mathfrak e}} 
\usepackage{scalerel}

\newcommand{\hatzero}{\widehat{0}}
\newcommand{\p}{\mathfrak{p}}

\NewDocumentCommand{\ssub}{O{0pt} O{.8} m t! e{_^}}{
  #3%
  \IfValueT{#5}{
    \IfBooleanTF{#4}{\sb{\hspace{#1}\scaleobj{#2}{#5}}}{\sb{#5}}
  }
  \IfValueT{#6}{\sp{#6}}
}
\RenewDocumentCommand{\ss}{O{0pt} O{0pt} O{.8} m e{_^}}{
  #4%
  \IfValueT{#5}{
    \sb{\hspace{#1}\scaleobj{#3}{#5}}
  }
  \IfValueT{#6}{
    \sp{\hspace{#2}\scaleobj{#3}{#6}}
  }
}

\newcommand{\cH}{\mathcal H}

\newcommand{\cY}{\mathcal Y}
\renewcommand{\setminus}{\smallsetminus}

\begin{document}


\title{Tropical function fields, finite generation, and faithful tropicalization}
\author{Omid Amini}
\address{CNRS - CMLS, \'Ecole Polytechnique, Palaiseau, France}
\email{omid.amini@polytechnique.edu}

\author{Shu Kawaguchi}
\address{Department of Mathematics, Graduate School of Science, Kyoto University, Kyoto 606-8502, Japan}
\email{kawaguch@math.kyoto-u.ac.jp}

\author{JuAe Song}
\address{Department of Mathematics, Graduate School of Science, Kyoto University, Kyoto 606-8502, Japan}
\email{song.juae.8m@kyoto-u.ac.jp}


\begin{abstract} Given an algebraic variety defined over a discrete valuation field and a skeleton of its Berkovich analytification, the tropicalization process transforms function field of the variety to a semifield of tropical functions on the skeleton. Our main result offers a purely polyhedral characterization of this semifield: we show that a tropical function is in the image of the tropicalization map if and only if it takes the same slope near infinity along parallel half-lines of the skeleton. This extends a result of Baker and Rabinoff in dimension one to arbitrary dimensions. We use this characterization to establish that this semifield is finitely generated over the semifield of tropical rational numbers, providing a new proof of a recent result by Ducros, Hrushovski, Loeser and Ye in the discrete valued field case.
As a second application, we present a new proof of the faithful tropicalization theorem by Gubler, Rabinoff and Werner in the discrete valuation field case. The proof is constructive and provides explicit coordinate functions for the embedding of the skeleton, extending the existing results in dimension one to arbitrary dimensions. 
\end{abstract}

\maketitle

\tableofcontents

\section{Introduction} Throughout this paper, let $K$ be a complete discrete valuation field with valuation $\ss \nu_K\colon K \to \ZZ\cup\{\infty\}$, valuation ring $R$ with maximal ideal $\mathfrak{m}$, and residue field $\kappa$. 

Let $X$ be a smooth projective geometrically connected variety over $K$ and let $\Xscr$ be a strictly semistable model of $X$, defined over the valuation ring of some finite extension of $K$. Let $\cH$ denote a horizontal divisor of $\Xscr$ such that $(\Xscr, \cH)$ forms a strictly semistable pair. The special fiber of $\Xscr$ is denoted by $\Xscr_\kappa$. After performing a possible finite base change, we assume that each connected component of the intersection of any collection of irreducible components of $\cH$ is geometrically connected. 

Let $\Delta(\Xscr, \cH)$ be the dual complex of $\Xscr_\kappa + \cH$, an abstract polyhedral complex with a recession fan defined as follows. 
Let $Z_v$ for $v \in V_\f$ be the irreducible components of $\Xscr_\kappa$ and $W_\rho$ for $\rho \in \ss V_\infty$ the irreducible components of $\cH$. To any connected component of an intersection of the form 
\[
\bigcap_{v \in A} Z_v \cap \bigcap_{\rho \in B} W_\rho
\]
where $A \subseteq V_\f$ is nonempty and  $B \subseteq \ss V_\infty$, we associate a copy of the polyhedron $\sigma^A \times \RR_{\geq 0}^{B}$. Here, $\sigma^A$ is the standard simplex in $\RR^{A}$ consisting of points with non-negative coordinates summing to one. The collection of these polyhedra forms $\Delta(\Xscr, \cH)$, which glue together to define the polyhedral space $\supp{\Delta(\Xscr, \cH)}$ to which we refer as the support of $\Delta(\Xscr, \cH)$. Moreover, the collection of cones consisting of a copy of $\RR^B_{\geq 0}$ associated with each connected component of a nonempty intersection of the form $\bigcap_{\rho \in B} W_\rho$, as above, forms an abstract fan $\Upsilon$ referred to as the \emph{recession fan} of $\Delta(\Xscr, \cH)$. These cones glue together and form the support $\supp{\Upsilon}$ of  $\Upsilon$. In particular, the rays of $\Upsilon$ correspond to the components $W_\rho$ of $\cH$. See Sections~\ref{subsubsec:abstract:polyherdal:complex} and~\ref{sec:tropical} for more details. 

Two rays in $\Delta(\Xscr, \cH)$ are said to be \emph{parallel} if they are associated with the same component of $\cH$, meaning they correspond to the same ray in $\Upsilon$. More generally, two half-lines in $\supp{\Delta(\Xscr, \cH)}$ are \emph{parallel} if the corresponding half-lines issued from the origin in $\supp{\Upsilon}$ are the same. 

Let $X^{\an}$ denote the Berkovich analytification of $X$ (see \cite{Ber90}). The dual complex $\Delta(\Xscr, \cH)$ (more precisely, its support) naturally embeds in $X^{\an}$ forming a subspace called the \emph{skeleton $\Sk(\Xscr, \cH)$ of $X^{\an}$ associated with $(\Xscr, \cH)$}. This allows to identify $\supp{\Delta(\Xscr, \cH)}$ with $\Sk(\Xscr, \cH) \hookrightarrow X^{\an}$.

Let $\Rat(X)$ denote the function field of $X$.  Identifying $\supp{\Delta(\Xscr, \cH)}$ with the skeleton $\Sk(\Xscr, \cH) \hookrightarrow X^{\an}$, we associate to each point $x \in \supp{\Delta(\Xscr, \cH)}$ a valuation 
\[
\nu_x\colon \Rat(X) \to \RR\cup\{\infty\}.
\]
For any $f \in \Rat(X)$, the (\emph{abstract}) \emph{tropicalization} of $f$ is defined as 
\[
\trop(f)\colon \supp{\Delta(\Xscr, \cH)} \to \RR \cup\{\infty\}, \qquad x \mapsto \nu_x(f).
\]   
More generally, consider the base change $X_{\bar{K}}$  of $X$ to an algebraic closure $\bar{K}$ of $K$. The tropicalization map naturally extends to $\Rat(X_{\bar K})$, and for $f\in \Rat(X_{\bar K})$, it yields a function $\trop(f)$, as described above, on $\supp{\Delta(\Xscr, \cH)}$.

The set $\trop(\Rat(X_{\bar K}))$ of tropicalized rational functions on $X_{\bar K}$ is a semifield with the operation of addition corresponding to taking pointwise minimum of functions, and the operation of multiplication corresponding to the usual addition.

Our first main result gives a characterization of this semifield. This is done as follows. 

We define $\Rat(\Delta(\Xscr, \cH))$ as the union of $\{\infty\}$ and the space of facewise integral $\QQ$-affine functions on $\Delta(\Xscr, \cH)$ which take the same slope near infinity along parallel rays of $\Delta(\Xscr, \cH)$.  
Similarly, we 
define $\ss\Rat(\supp{\Delta(\Xscr, \cH)})$ as the union of $\{\infty\}$ and the space of piecewise integral $\QQ$-affine functions on $\Delta(\Xscr, \cH)$ which exhibit the same slope near infinity along parallel half-lines of $\supp{\Delta(\Xscr, \cH)}$. 
Equivalently, $\ss\Rat(\supp{\Delta(\Xscr, \cH)})$ is the union of $\Rat(\Delta^\prime)$ over integral $\QQ$-affine refinements $\Delta^\prime$ of $\Delta(\Xscr, \cH)$ (see Section~\ref{subsubsec:abstract:polyherdal:complex}). We refer to $\ss\Rat(\supp{\Delta(\Xscr, \cH)})$ as the {\em tropical function field of $(\Xscr, \cH)$}.

\begin{Theorem}[Characterization Theorem]
\label{thm:Q:1}
For every $f\in \Rat(X_{\bar K})$, the tropicalization $\trop(f)$ belongs to $\ss\Rat(\supp{\Delta(\Xscr, \cH)})$. The resulting map
\begin{equation}
\label{eqn:trop:K(X):to:M}
\trop\colon \Rat(X_{\bar{K}}) \to \ss\Rat(\supp{\Delta(\Xscr, \cH)})
\end{equation}
is surjective. 

Consequently, the tropical function field coincides with the semifield of tropicalized rational functions on $X_{\bar K}$. 
\end{Theorem}

This result generalizes the work of Baker and Rabinoff  \cite{BR15} in dimension one to arbitrary dimensions, while also accounting for the presence of horizontal divisors.  See also the work of Amini and Iriarte \cite{AI22}, which addresses the case of trivially valued fields of characteristic zero, and establishes the surjectivity of an analogous tropicalization map, defined relative to a compactification.  
We note that our proof of the characterization theorem gives the surjectivity 
of tropicalization map from $\Rat(X_{L})$ to $\ss\Rat(\supp{\Delta(\Xscr, \cH)})$ for an extension $L$ of $K$ 
with value group $\QQ$. 

\smallskip

  To state our second main result, we recall that a semifield $\mathbb F$ is said to be finitely generated over a semifield $\mathbb E$ if there exist finitely many elements $a_1, a_2, \ldots, a_n$ of $\mathbb F$ such that every element of $\mathbb F$ can be expressed as  $p(a_1, \ldots, a_n)/q(a_1, \ldots, a_n)$ with $p$ and $q$  $n$-variable polynomials over $\mathbb E$. 

\smallskip

Let $\TT\QQ\colonequals \QQ \cup \{\infty\}$ denote the semifield of rational numbers with the operations of min for addition and $+$ for multiplication.  We observe that the definition of tropical function fields extends to any abstract polyhedral complex endowed with a recession fan (see Section~\ref{subsubsec:abstract:polyherdal:complex}).

\begin{Theorem}[Finite generation of tropical function fields]
\label{thm:Q:2}
Let $\Delta$ be an abstract polyhedral complex with a  recession fan. The tropical function field
  $\ss\Rat(\supp{\Delta})$ is finitely generated over $\TT\QQ$. 
\end{Theorem}

In particular, $\ss\Rat(\supp{\Delta(\Xscr, \cH)})$ as above is finitely generated over $\TT\QQ$. In a recent work \cite{DHLY24}, using model theory of valued fields, Ducros, Hrushovski, Loeser, and Ye prove that $\trop(\Rat(X_{\bar K}))$ is finitely generated as a semifield over the value group. A combination of our Theorems~\ref{thm:Q:1} and~\ref{thm:Q:2} yields a new proof of their result in the discrete valued field case,
showing in addition the finite generation of $\trop(\Rat(X_{L}))$ for an extension $L$ of $K$ with value group $\QQ$ without requiring $L$ to be defectless. Note that as it was stated in \cite{DHLY24}, finite generation property does not hold in general for $\trop(\Rat(X))$. 

\smallskip
 
As an application of Theorems~\ref{thm:Q:1} and  \ref{thm:Q:2}, we provide an alternative proof of the faithful tropicalization theorem of Gubler, Rabinoff, and Werner \cite{GRW16} in the discrete valued  field case. See also~\cite[Thm~5.1]{Duc12} and~\cite[\S1.2]{DHLY24}.

\begin{Theorem}[{Gubler--Rabinoff--Werner \cite{GRW16}}]
\label{thm:GRW}
Let $X$ be a smooth projective variety over $K$ with a strictly semistable pair $(\Xscr, \cH)$, as above. Then, there exist $f_1, f_2, \ldots, f_n \in \Rat(X_{\bar{K}})$ such that, 
setting $U \colonequals X_{\bar{K}} \setminus \bigcup_{i=1}^n \supp{\zero(f_i)}$, the map 
\[
\psi\colon U^{\an} \to \RR^n, \quad 
x \mapsto (\nu_x(f_1), \ldots, \nu_x(f_n))
\]
provides a faithful tropicalization of the skeleton $\Sk(\Xscr, \cH)$ of $(\Xscr, \cH)$, that is, the restriction of $\psi$ to $\Sk(\Xscr, \cH)$ is a homeomorphism onto its image that preserves the integral structures.  
\end{Theorem}

Our proofs of Theorems~\ref{thm:Q:1},~\ref{thm:Q:2} and~\ref{thm:GRW} yield additional insights. We explicitly construct $F_1, \ldots, F_n \in \ss\Rat(\supp{\Delta(\Xscr, \cH)})$ such that 
the map 
\[
(F_1, \ldots, F_n)\colon \supp{\Delta(\Xscr, \cH)} \to \RR^n
\]
is a homeomorphism onto its image preserving the integral structures.  By Theorem~\ref{thm:Q:1}, there exist $f_1, \ldots, f_n \in \Rat(X_{\bar{K}})$ satisfying $\trop(f_1) = F_1, \ldots, \trop(f_n) = F_n$. 
Our proof shows that the embedding given by
$(f_1, \ldots, f_n)$ provides a faithful tropicalization of the skeleton $\Sk(\Xscr, \cH)$ as stated in Theorem~\ref{thm:GRW}. Moreover, $\trop(f_1), \ldots, \trop(f_n)$ generate 
$\trop(\Rat(X_{\bar K}))$, which is equal to $\ss\Rat(\supp{\Delta(\Xscr, \cH)})$ by Theorem~\ref{thm:Q:1}.  In particular, the proof provides explicit coordinate functions for the embedding of the skeleton. This way, we extend the existing results on effective tropicalization in dimension one~\cite{BPR16, CM19, Helm19, KMM08, KY21,  MRS23}, and for special varieties of higher dimension~\cite{CHW14, DP16, Kut16, KY19, Ran17}, to arbitrary varieties in all dimension.

\smallskip

Our proof of Theorem~\ref{thm:Q:2} gives a stronger result. Let $\Sigma$ be an integral $\QQ$-affine polyhedral complex in $\RR^n$. 
We define $\ss \Rat(\supp{\Sigma})$ to be the set consisting of $\infty$ and piecewise integral $\QQ$-affine functions on $\supp{\Sigma}$ which that take the same slope near infinity along any pair of parallel half-lines in $\supp{\Sigma}$ (see Definition~\ref{def:rat:Delta}).

Let $\Rat(\RR^n)$ denote the subsemifield of the semifield of piecewise integral $\QQ$-affine functions on $\RR^n$ generated by the coordinate functions $\x_1, \dots, \x_n$.
Then, we have a natural restriction map 
\[
\pi\colon \Rat(\RR^n) \to \ss\Rat(\supp{\Sigma}).
\]
We denote by ${\rm Ker}(\pi)$ the kernel of $\pi$ as a congruence. For the notion of a congruence, which serves as the analogue of an ideal when forming a quotient in the context of semirings, see Section~\ref{subsec:congruences}. 

\begin{Theorem}
\label{thm:Q:2:stronger}
The restriction map $\pi$ is surjective, meaning that the restrictions of the coordinate functions $\x_1, \ldots, \x_n$ to
$\supp{\Sigma}$ generate the semifield $\ss\Rat(\supp{\Sigma})$ over $\TT\QQ$. Moreover,  ${\rm Ker}(\pi)$ is finitely generated as a congruence. 
\end{Theorem}

This generalizes a result of Song \cite{Song} in dimension one
to arbitrary dimensions
when $\Sigma$ is defined over $\TT\QQ$.

\medskip
As a final remark, although we have restricted our focus to the case of discretely valued fields in order to avoid technicalities and have a clearer presentation, we are confident that the arguments employed in this work can be extended to more general valued fields. Moreover, the condition on $X$ being smooth,  projective and geometrically connected is not restrictive, since our results on tropicalizations concern birational geometry of varieties over an algebraic closure of $K$. 

\medskip

The text is organized as follows. In Section~\ref{subsec:polyhedral:complexes}, we provide preliminaries on polyhedral geometry, followed by background materials on tropical and non-Archimedean geometry in Section~\ref{sec:tropical}. Theorem~\ref{thm:Q:1} is proved in Section~\ref{sec:proof:Q1}. The proofs of Theorems~\ref{thm:Q:2},~\ref{thm:GRW}, and~\ref{thm:Q:2:stronger} are given in Section~\ref{sec:strategy}.

\bigskip
{\sl Acknowledgment.}\quad
This project started during the workshop ``Nonarchimedean geometry and related fields'' in Kyoto in January 2024. We thank the organizers of the workshop and Kyoto University for 
hospitality. We thank Matthieu Piquerez and Kazuhiko Yamaki for helpful discussions.
O.A. is part of the ANR project AdAnAr (ANR-24-CE40-6184). S.K. is partially supported by KAKENHI 23K03041 and 20H00111.


\section{Polyhedral complexes}
\label{subsec:polyhedral:complexes}
We gather some definition and notation regarding embedded and abstract polyhedral complexes. 
\subsection{Polyhedral complexes in $\RR^n$}
We denote the points of $\RR^n$ with bold letters $\mathbf{u}$ with coordinates $u_1, \dots, u_n$. 

Let $\langle\cdot\,, \cdot \rangle$ denote the standard inner product of $\RR^n$. An {\em integral $\QQ$-affine function} on $\RR^n$ is a function 
of the form 
\[
\RR^n \to \RR, \quad 
\mathbf{u} \mapsto \langle \mathbf{m}, \mathbf{u}\rangle+\gamma
\]
for some $\mathbf{m} \in \ZZ^n$ and $\gamma \in \QQ$. A function of the above form with $\mathbf{m} \in \RR^n$ and $\gamma\in \RR$ is called an affine function.

An {\em integral $\QQ$-affine polyhedron} in $\RR^n$ is a subset of the form 
\[
\sigma = 
\bigl\{
\mathbf{u} \in \RR^n \mid L_i(\mathbf{u}) \geq 0 \; \text{for $i = 1, \ldots ,r$}
\bigr\}
\]
for integral $\QQ$-affine functions $L_i$ on $\RR^n$.  Each face of an integral $\QQ$-affine polyhedron is again integral $\QQ$-affine. When $\sigma$ is compact, it is called an \emph{integral $\QQ$-affine polytope}. If the $L_i$s in the definition of $\sigma$ are integral linear, then $\sigma$ is called a \emph{rational cone}.

An {\em integral $\QQ$-affine polyhedral complex} $\Sigma$
in $\RR^n$ is a polyhedral complex whose cells are integral $\QQ$-affine. Recall that this means $\Sigma$ is a finite collection of integral $\QQ$-affine polyhedra $\sigma$ such that 
\begin{itemize}
    \item for each $\sigma \in \Sigma$, any face $\tau$ of $\sigma$ is included in $\Sigma$.
    \item for each pair $\sigma,\eta \in \Sigma$ with nonempty intersection, $\sigma\cap \eta$ is a common face of $\sigma$ and $\eta$.
\end{itemize}

Let $\Sigma$ be an integral $\QQ$-affine polyhedral complex in $\RR^n$. For $\sigma, \tau \in \Sigma$, we write $\tau \preceq \sigma$ if $\tau$ is a face of $\sigma$.  A \emph{ray} of $\Sigma$ is an unbounded one-dimensional element $\varrho$ in $\Sigma$. The {\em support} of $\Sigma$ denoted by $\supp{\Sigma}$ is the union of polyhedra $\sigma\in \Sigma$. We say that $\Sigma$ is \emph{complete} if $\supp{\Sigma}=\RR^n$.

Each integral $\QQ$-affine polyhedron $\sigma$ in $\RR^n$ can be expressed as a sum $P+\tau$, where $P$ is an integral $\QQ$-affine polytope and $\tau$ is a rational cone, uniquely determined by $\sigma$. We refer to $\tau$ as the \emph{recession fan} of $\sigma$, and denote it by $\rec(\sigma)$. When the collection $\rec(\sigma)$, for $\sigma \in \Sigma$, forms a fan, we call it the \emph{recession fan} of $\Sigma$, see~\cite{BS11}.

Let $\Sigma, \Sigma^\prime$ be two integral $\QQ$-affine polyhedral complexes in $\RR^n$. We say that $\Sigma^\prime$ is a {\em refinement} of 
$\Sigma$ if $\supp{\Sigma} =\supp{\Sigma^\prime}$ and for each $\sigma^\prime \in \Sigma^\prime$, there exists 
$\sigma \in \Sigma$ with $\sigma^\prime \subseteq \sigma$.

Let $\sigma \subseteq \RR^n$ be an integral $\QQ$-affine polyhedron. An ({\em integral $\QQ$-}){\em affine function on $\sigma$} is by definition the restriction to $\sigma$ of an (integral $\QQ$-)affine function $\RR^n \to \RR$. 

\smallskip

We define the space of facewise, resp. piecewise,  (integral $\QQ$-)affine functions on $\Sigma$, resp. $\supp{\Sigma}$, as follows. 

\begin{Definition}[Facewise (integral $\QQ$-)affine functions and $\Rat(\Sigma)$]
    \label{def:facewise}
    Let $\Sigma$ be an integral $\QQ$-affine polyhedral complex in $\RR^n$. 
    We say that $F \colon \supp\Sigma\to\RR$ is {\em facewise} (integral $\QQ$-)affine on $\Sigma$ if $\rest{F}{\sigma}\colon \sigma \to \RR$ is (integral $\QQ$-)affine for any $\sigma \in \Sigma$. 

   We denote by $\Rat(\Sigma)$ the union of $\{\infty\}$ and the space of facewise integral $\QQ$-affine functions on $\Sigma$ which take the same slope at infinity along any pair of parallel rays in $\Sigma$. Similarly, we define $\Rat(\Sigma,\RR)$ as the union of $\{\infty\}$ and the space of all facewise affine functions on $\Sigma$ which take the same slope at infinity along any pair of parallel rays in $\Sigma$.
    \end{Definition}

\begin{Definition}[Piecewise (integral $\QQ$-)affine functions and $\ss\Rat(\supp\Sigma)$]
\label{def:rat:Delta}
Let $\Sigma$ be an integral $\QQ$-affine polyhedral complex in $\RR^n$.  
We say that $F \colon \supp\Sigma\to\RR$ is \emph{piecewise} integral $\QQ$-affine if there exists an integral $\QQ$-affine refinement $\Sigma^\prime$ of $\Sigma$ such that the restriction
$\rest{F}{\sigma^\prime}\colon \sigma^\prime \to \RR$ is an integral $\QQ$-affine function on $\sigma^\prime$ for each $\sigma^\prime \in \Sigma^\prime$.
We denote by $\ss\Rat(\supp\Sigma)$ the union of $\{\infty\}$ and the space of all piecewise integral $\QQ$-affine functions $F$, as above, where $F$ takes the same slopes at infinity along any two parallel half-lines of $\supp{\Sigma}$.  

Similarly, we define \emph{piecewise} affine functions on $\supp{\Sigma}$ and $\ss \Rat(\supp\Sigma, \RR)$, by removing the integrality constraint. \qedhere
\end{Definition}

Note that each element $F \in \Rat(\Sigma)$ belongs to $\Rat(\Sigma^\prime)$ for some integral $\QQ$-affine refinement of $\Sigma^\prime$. Moreover, a pair of integral $\QQ$-affine refinements $\Sigma^{\prime}$ and $\Sigma^{\dprime}$ have a common integral $\QQ$-affine refinement.  Therefore, $\Rat(\supp\Sigma)$ is closed under addition and taking minimum. This turns $\Rat(\supp\Sigma)$ into a semifield (see Section~\ref{subsec:tropical:semifield}).  The same remark applies to $\ss \Rat(\supp{\Sigma},\RR)$.

\smallskip

For future use, we define convex functions on integral $\QQ$-affine polyhedral complexes.  

\begin{Definition}[Convex function]
\label{def:convex:function}
Let $\Sigma$ be an integral $\QQ$-affine polyhedral complex in $\RR^n$ and let $H \colon \supp{\Sigma}\to\RR$ be a facewise (integral $\QQ$-)affine function on $\Sigma$. We say that $H$ is {\em convex} if, for any $\sigma \in \Sigma$, there exists an (integral $\QQ$-)affine function $L_\sigma\colon \RR^n \to \RR$ such that $\rest{H}{\sigma} = \rest{L_\sigma}{\sigma}$ and   
\[
\rest{(H-L_\sigma)}{\eta\setminus \sigma} \geq 0
\]
on $\eta\setminus\sigma$ for each $\eta \in \Sigma$ with $\eta \succeq\sigma$. We say that $H$ is \emph{strictly convex} if the inequalities above are strict on $\eta \setminus \sigma$.
\end{Definition}

We say that a function $H \colon \supp{\Sigma}\to\RR$
is {\em concave} if $-H$ is convex.

\begin{Definition}[Convexity of a polyhedral complex]
\label{def:convex:polyhedral:complex}
An integral $\QQ$-affine polyhedral complex $\Sigma$ is said to be {\em convex} if 
there exists a strictly convex function $H \in \Rat(\Sigma)$.
\end{Definition}

\subsection{Abstract polyhedral complex with a recession fan}
\label{subsubsec:abstract:polyherdal:complex}
For a set $V$, we denote by $\ssub{2}!^V$ the power set of $V$.

An \emph{abstract polyhedral complex with a recession fan} is defined as follows. Let $V_\f$ and $\ss V_\infty$ be two disjoint finite sets, with $V_\f$ nonempty. Let $\Pi$ be a finite delta-complex on the set $V_\f \cup \ss V_\infty$. This means that $\Pi$ is a finite set endowed with a partial order $\preceq$ and there exists a map $\zeta\colon \Pi \to \ssub{2}!^{V_\f\cup \ss V_\infty}$ with the following properties: 
\begin{enumerate}[label=($\arabic*$)]
\item\label{Pi1} There exists a minimum element in $\Pi$, denoted by $\hatzero$, with $\zeta(\hatzero)=\emptyset$.
\item \label{Pi2} For each $x \in V_\f \cup \ss V_\infty$, there exists a unique element $\tau$ of $\Pi$ with $\zeta(\tau)=\{x\}$. 
\item\label{Pi3} For each $\tau \in \Pi$, the restriction of $\zeta$ to the interval $[\hatzero, \tau] \subseteq \Pi$ is a bijection with the poset of subsets of $\zeta(\tau)$. Here, the partial order on subsets is given by inclusion. 
\end{enumerate}
In particular, for each pair of elements $\tau \preceq \eta$ in $\Pi$, we have $\zeta(\tau) \subseteq \zeta(\eta)$.

Denote by $\Delta$ the set of all elements of $\Pi$ with nonempty intersection with $V_\f$. We call $\Delta$ an \emph{abstract polyhedral complex with vertex set $V_\f$ and recession rays $\ss V_\infty$}. By an abuse of the notation, we denote the restriction of $\zeta$ to $\Delta$ by $\zeta \colon \Delta \to \ssub{2}!^{V_\f\cup \ss V_\infty}$. The complement $\Pi \setminus \Delta$ is itself a delta-complex. Let $\Upsilon$ be the abstract polyhedral complex with vertex set $\{0\}$ and recession rays $\ss V_\infty$ which as a set is equal to $\Pi \setminus \Delta$ endowed with the same partial order $\preceq$, with the map $\ss \zeta_\infty \colon \Upsilon \to \ssub{2}!^{\{0\} \cup \ss V_\infty}$ defined by $\ss\zeta_\infty(\tau) = \{0\}\cup \zeta(\tau)$ for all $\tau \in \Pi \setminus \Delta$. Although the definition of $\Upsilon$ is dependent on $\Pi$, by an abuse of the notation we refer to $\Upsilon$ as the recession fan of $\Delta$, and call $\Delta$ \emph{an abstract polyhedral complex with recession fan $\Upsilon$}. 

\smallskip
Given a nonempty subset $A\subseteq V_\f$ and a subset $B \subseteq  \ss V_\infty$, we define the {\em extended simplex}
\[
\ss\Theta_{A\cup B} \colonequals \sigma^{A} \times \RR_{\geq 0}^{B} \subseteq \RR^{A\cup B} \subseteq \RR^{V_\f \cup \ss V_\infty},
\]
where $\RR_{\geq 0}^{B}$ denotes the non-negative orthant in $\RR^B$, and $\sigma^A$ denotes the standard simplex in $\RR^A$,
\begin{equation}
\label{eqn:standard:simplex}
\sigma^{A} \colonequals 
\Bigl\{
\ss{(u_v)}_{v\in A} \in \RR^{A}_{\geq 0} \, \mid \,
\sum_{v\in A} u_v = 1
\Bigr\}.
\end{equation}

 We denote by $\ss\e_\rho$, $\rho\in \ss V_\infty$, the standard basis of $\RR^{\ss V_\infty}$. 

Given an abstract polyhedral complex $\Delta$ with recession fan $\Upsilon$ as defined above, we associate to each element $\tau \in \Delta$ a copy $\Theta_\tau$ of the polyhedron $\Theta_{\zeta(\tau)} \subset \RR^{V_\f \cup \ss V_\infty}$. The property~\ref{Pi3} ensures that the polyhedra $\Theta_{\tau}$ are nonempty and they glue together, along inclusions, for each pair of elements $\tau \preceq \eta$ in $\Delta$, of $\Theta_{\tau}$ as a face of $\Theta_{\eta}$. 

The resulting space is denoted by $\supp{\Delta}$ and is called the \emph{support} of $\Delta$. When there is no risk of confusion, we simplify the notation and write $\tau$ in place of $\Theta_{\tau}$. 

Using this simplified notation:
\begin{itemize}
    \item The \emph{finite part} of $\tau$ is denoted by $\tau_\f$ and corresponds to the simplex $\sigma^{\zeta(\tau)\cap V_\f}$.
    \item The \emph{infinite part} of $\tau$ is denoted by $\ss\tau_\infty$, it corresponds to the cone $\RR_{\geq 0}^{\zeta(\tau)\cap \ss V_\infty}$ and naturally belongs to $\Upsilon$.
\end{itemize} 
Similarly, we define $\supp{\Upsilon}$, the support of the recession fan $\Upsilon$,  by gluing the collection of $\Theta_\tau$, identified with $\RR_{\geq 0}^{\ss \zeta_\infty(\tau)\cap \ss V_\infty}$, for all $\tau\in \Upsilon$. (This can be viewed as a multifan in $\RR^{\ss V_{\infty}}$.)

A \emph{ray} in $\Delta$ is an element $\varrho$ with $\zeta(\varrho)\cap V_\f$ and $\zeta(\varrho)\cap \ss V_\infty$ each have  size one. Two rays $\varrho$ and $\varrho^\prime$ are {\em parallel} if they correspond to the same element of $\ss V_\infty$.

More generally, two half-lines in $\supp{\Delta}$ are 
\emph{parallel} if they correspond to the same half-line issued from the origin in $\supp{\Upsilon}$. 
 
\begin{Example}
\begin{enumerate}
    \item[(1)]
Let $V_\f = \{v_1, v_2\}$ and $\ss V_\infty = \{\rho, \rho^\prime\}$. We set 
\[
\Delta = \{v_1, v_2, e, \varrho_1, \varrho_1^\prime, \varrho_2, \varrho_2^\prime, \sigma_1, \sigma_2\}, \qquad 
\Upsilon = \{ \widehat{0}, \rho, \rho^\prime, \tau_{1}, \tau_{2}\}, \qquad
\Pi = \Delta \cup \Upsilon
\]
and set 
\begin{align*}
& \zeta(v_i)  = \{v_i\}, \quad  \zeta(e) = \{v_1, v_2\}, \quad  
\zeta(\varrho_i) = \{v_i, \rho\}, \quad  
\zeta(\varrho_i^\prime) = \{v_i, \rho^\prime\}, \quad  
\zeta(\sigma_i) = \{v_i, \rho, \rho^\prime\}, 
\\
& 
\zeta(\rho) = \{\rho\}, \quad  \zeta(\rho^\prime) = \{\rho^\prime\}, \quad 
\zeta(\tau_{1})  = \zeta(\tau_{2})=\{\rho, \rho^\prime\}. 
\end{align*}
Then, a half-line in $\sigma_1$ is parallel to a half-line in $\sigma_2$ if they are either parallel to $\rho$ or they are parallel to $\rho^\prime$. 

    \item[(2)]
We take the same example as above, but we change $\Upsilon$ to 
\[
\Upsilon = \{ \widehat{0}, \rho, \rho^\prime, \tau\}
\quad\text{and}\quad \zeta(\tau) = \{\rho, \rho^\prime\}. 
\]
In this second example, the recession cones of $\sigma_1$ and $\sigma_2$ are the same, and so 
a half-line in $\sigma_1$ is parallel to a half-line in $\sigma_2$ if they are the same half-line in $\tau$. \qedhere
\end{enumerate}
\end{Example}

An \emph{integral $\QQ$-affine refinement} $\Delta^\prime$ of $\Delta$ consists in a collection of integral $\QQ$-affine refinements of each $\tau$, in the sense of the previous section when $\tau$ is viewed as $\Theta_{\zeta(\tau)}$ embedded in $\RR^{V_\f\cup \ss V_\infty}$. These refinements are required to be compatible with each other along inclusion of faces. 
Two rays in $\Delta^\prime$ are parallel 
if they are parallel half-lines in $\supp{\Delta} = \supp{\Delta^\prime}$.

Given an abstract polyhedral complex $\Delta$ with  recession fan $\Upsilon$, we define $\Rat(\Delta)$ as the union of $\{\infty\}$ and the space of facewise integral $\QQ$-affine functions on $\Delta$ which take the same slope near infinity along parallel rays of $\Delta$. We define $\ss\Rat(\supp{\Delta})$ as the union of $\{\infty\}$ and the collection of all facewise integral $\QQ$-affine functions on some integral $\QQ$-affine refinement $\Delta^\prime$ which take the same slopes near infinity along parallel rays of $\Delta^\prime$. Equivalently, $\ss\Rat(\supp{\Delta})$ is the union of 
$\{\infty\}$ and the collection of all piecewise integral  $\QQ$-affine functions on $\supp{\Delta}$ which take the same slope near infinity along two parallel half-lines in $\supp{\Delta}$. The set $\ss\Rat(\supp{\Delta})$ has the structure of a semifield with addition corresponding to taking pointwise minimum of functions, and multiplication corresponding to the usual addition (see Section~\ref{subsec:tropical:semifield}).

\section{Tropical preliminaries} \label{sec:tropical}
The purpose of this section is primarily to provide definitions and establish notation. For further details, we refer the reader to \cite{Ber90, Ber99, Ber04, GRW16, BFJ16, dJ, KY19, MN15}. 

\subsection{Notation and conventions}
\label{subsec:notation}
As in the introduction, let $K$ be a complete discrete valuation field. Let $R$ be the valuation ring of $K$, $\mathfrak{m}$ the maximal ideal of $R$, and $\kappa = R/\mathfrak{m}$ the residue field. 
We fix a uniformizer $\varpi \in \mathfrak{m}$ of $R$. Let $\ss \nu_K\colon K^\times \to \ZZ$ be the valuation map with $\ss \nu_K(\varpi) = 1$. 

By a {\em variety $X$ over $K$} we mean an irreducible, reduced, and separated scheme $X$ of finite type over $K$. By a {\em variety $\Xscr$ over $R$}, we mean an irreducible, reduced, and separated scheme $\Xscr$, flat and of finite type over $R$. The special fiber of $\Xscr$ is denoted by $\Xscr_\kappa \colonequals \Xscr \times_{R} \kappa$, and the generic fiber of $\Xscr$ is denoted by $\Xscr_K \colonequals \Xscr \times_R K$.

\subsection{Semistable models and skeleta}
\label{subsec:semistable:models}

Let $\Xscr$ be a projective variety over $R$ of dimension $d+1$. Let $Z_{v}$, $v\in V_\f$, denote the irreducible components of the special fiber $\Xscr_\kappa$. For a subset $A \subseteq V_\f$, we put $Z_{A} = \bigcap_{v \in A} Z_{v}$. We say that 
$\Xscr$ is  {\em strictly semistable} if $\Xscr$ is covered by open subsets $\Ucal$ that admit an \'etale morphism 
\[
\psi\colon \Ucal \to \Sscr \colonequals \Spec(R[\x_0, \ldots, \x_d]/\langle \x_0 \cdots \x_r -\varpi\rangle) \quad \text{for } 0 \leq r \leq d.
\]
Under these conditions, $\Xscr$ is regular, the special fiber $\Xscr_\kappa$ is reduced, and for each subset $A \subseteq V_\f$, each connected component of $Z_{A}$ is smooth over $\kappa$. For further details, see \cite[Definition~3,1]{GRW16} and \cite[2.16]{dJ}.  

Let $X$ be a smooth projective variety of dimension $d$ over $K$. A {\em strictly semistable model} of $X$ over $R$ is a strictly semistable projective variety $\Xscr$ over $R$ endowed with an isomorphism $\Xscr_K \to X$ over $K$.

Let $\Xscr$ be a projective variety over $R$ of dimension $d+1$. Let $\cH = \sum_{\rho\in \ss V_\infty} W_\rho$ be a sum of finitely many hypersurfaces $W_\rho$ on $\Xscr$, where the $W_\rho$ are irreducible and pairwise distinct.  We set $\cY \colonequals  \Xscr_\kappa + \cH$.

We say that the pair $(\Xscr, \cH)$ is a {\em strictly semistable pair} if 
$\Xscr$ is covered by open subsets $\Ucal$ which admit an \'etale morphism 
\[
\psi\colon 
\Ucal \to \Sscr \colonequals \Spec(R[\x_0, \ldots, \x_d]/\langle \x_0 \cdots \x_r - \varpi\rangle)
\]
for $0 \leq r \leq d$ such that the restriction of $W_\rho$ to $\Ucal$ is either trivial or defined by $\psi^*(\x_j)$ for some $r+1 \leq j \leq d$.  In this case, $\Xscr$ is strictly semistable, $\cY$ is a divisor with strict normal crossings on $\Xscr$, and for each $B \subseteq \ss V_\infty$, $W_{B} = \bigcap_{\rho \in B} W_\rho$ is a disjoint union of varieties over $R$, which are themselves strictly semistable over $R$, cf. \cite[Definition~3,1]{GRW16} and \cite[6.3, 6.4]{dJ}. We refer to the $W_\rho$, $\rho\in \ss V_\infty$, as the {\em horizontal components} and to the $Z_{v}$, $v\in V_\f$, as the {\em vertical components} of $\cY$.

Given a smooth projective variety $X$ defined over $K$, we say that $(\Xscr, \cH)$ is a \emph{strictly semistable pair} of $X$ if $\Xscr$ is a strictly semistable model of $X$. 

From a strictly semistable pair $(\Xscr, \cH)$ of $X$, we obtain a delta-complex $\Pi(\Xscr, \cH)$ giving rise to an abstract polyhedral complex $\Delta(\Xscr, \cH)$ and a recession fan $\Upsilon(\Xscr, \cH)$, and the corresponding skeleton $\Sk(\Xscr, \cH)$ in $X^{\an}$, with its integral polyhedral complex structure,  as follows.

\subsubsection*{Stratification}
Let $(\Xscr, \cH)$ be a strictly semistable pair of $X$. 
By the definition of a strictly semistable pair,  
$\cY = {\Xscr}_\kappa + \cH$ is a reduced divisor on $\Xscr$, and each irreducible component of $\cY$ is either a connected smooth projective variety over $\kappa$ (for the $Z_{v}$, $v\in V_\f$), or a connected regular projective variety over $R$ (for the $W_{\rho}$, $\rho\in \ss V_\infty$). 

We define a stratification of $\cY$ as follows.
We put $\cY^{(0)} \colonequals \cY$, and for each $n \in \ZZ_{> 0}$, we let 
$\cY^{(n)}$ be the complement of the set of regular points in $\cY^{(n-1)}$. Thus, we obtain a chain of closed subsets: 
\[
\cY= \cY^{(0)}  \supsetneq 
\cY^{(1)}  \supsetneq 
\cdots \supsetneq 
\cY^{(t)}  \supsetneq 
\cY^{(t+1)}  = \emptyset, 
\]
where $t \leq d$. 

We say that $S$ is a {\em stratum} of $\cY$ if 
there exist subsets $A \subseteq V_\f$ and 
$B \subseteq \ss V_\infty$ such that 
$S$ is a nonempty connected component of $\left(Z_A \cap  Z_B\right) \setminus  \cY^{(|A|+|B|)}$. 
We say that $S$ is  {\em special} if $S  \subseteq \Xscr_\kappa$. 
Note that $S$ is a special stratum if and only if $S \neq \emptyset$ and $A \neq \emptyset$.

\subsubsection*{Dual complex and (extended) simplices associated with special strata} 
Let $\Pi(\Xscr, \cH)$ be the set consisting of an element $\widehat{0}$ and other elements in bijection with the strata of $\cY$. We endow this set with the partial order corresponding to the reverse inclusion of strata, and declaring that $\widehat{0}$ is the minimum element. We define the map $\zeta \colon \Pi(\Xscr, \cH) \to \ssub{2}^{V_\f\cup \ss V_\infty}$, by setting $\zeta(\widehat{0}) = \emptyset$ and sending each stratum to the union of the corresponding index sets $A \subseteq V_\f$ and $B\subseteq \ss V_\infty$. This turns $\Pi(\Xscr, \cH)$ into a delta-complex on the set $V_\f\cup \ss V_\infty$.

Let $\Delta(\Xscr, \cH)$ be the corresponding abstract polyhedral complex and $\Upsilon(\Xscr, \cH)$ the corresponding recession fan in the sense of Section~\ref{subsubsec:abstract:polyherdal:complex}. Specifically, $\Delta(\Xscr, \cH)$ is the set of elements $\tau$ in bijection with the special strata $S_\tau$ of $\cY$.  We refer to $\Delta(\Xscr, \cH)$ as the \emph{dual complex} of $(\Xscr, \cH)$.

Let $S_\tau$, $\tau \in \Delta(\Xscr, \cH)$, be a special stratum of $\cY$ with the corresponding index sets $A \subseteq V_\f$ and $B \subseteq \ss V_{\infty}$ as above. We associate to $\tau$ the extended simplex $\Theta_{\tau}=\Theta_{A\cup B}=\sigma^{A} \times \RR^{B}_{\geq 0}$.

There is a canonical embedding of $\Theta_{\tau}$ in $X^{\an}$, defined as follows. To each $\mathbf{u} \in \Theta_\tau$, we associate an absolute value $|\cdot|_{\tau, \mathbf{u}}$ on the function field $\Rat(X)$ which extends $|\cdot|_K$. Let $\xi = \xi_\tau$ be the generic point of $S_\tau$. For each $v \in A$ and $\rho\in B$, we choose a local equation $\z_v = 0$ of $Z_{v}$ and $\w_\rho = 0$ of $W_\rho$ in $\Xscr$ at $\xi$. 

Let $\widehat{\OO}_{\Xscr, \xi}$ be the completion of $\OO_{\Xscr, \xi}$. By \cite[\S2.4]{MN15}, each element $f \in \widehat{\OO}_{\Xscr, \xi}$ has an expansion of the form 
\[
f = \sum_{\mathbf{m} = (m_v)_{v\in A}\times (m_{\rho})_{\rho\in B} \in \ZZ^{A\cup B}_{\geq 0}} a_{\mathbf{m}} \prod_{v\in A}\z_{v}^{m_v} \cdot \prod_{\rho\in B}\w_{\rho}^{m_{\rho}}, 
\]
where each coefficient $a_{\mathbf{m}}$ is either zero or a unit in $\widehat{\OO}_{\Xscr, \xi}$. Such an expansion is called 
an {\em admissible expansion} of $f$. 

We write $\mathbf{u} = (u_v)_{v\in A}\times (u_\rho)_{\rho\in B}\in \sigma^{A} \times \RR^{B}_{\geq 0}$, and set 
\begin{equation}
\label{def:mathu:S}
  |f|_{\tau, \mathbf{u}} \colonequals \max_{\mathbf{m}  \in  \ZZ_{\geq 0}^{A\cup B}, \; a_{\mathbf{m}}  \neq 0} \Bigl\{ \exp\bigl(-\sum_{v\in A}u_{v} m_v -\sum_{\rho\in B} u_{\rho} m_\rho\bigr)
   \Bigr\}. 
\end{equation}
By \cite[Proposition~2.4.4]{MN15}, $|\cdot|_{\tau,\mathbf{u}}$ is a well-defined absolute value on $ \OO_{\Xscr, \xi}$. This absolute value extends to an absolute value $|\cdot|_{\tau, \mathbf{u}}$ on $\Rat(X)$. 
Furthermore, since $\varpi = \lambda \prod_{v\in A} \z_{v}$ for some unit $\lambda$ in $ \OO_{\Xscr, \xi}$, we have $|\varpi|_{\tau, \mathbf{u}} = \exp(-1)$,
which shows that
$|\cdot|_{\tau, \mathbf{u}}$ agrees with $|\cdot|_K$ for elements in $K$. 
Thus, for each $\mathbf{u} \in \Theta_\tau$, we have a point $|\cdot|_{\tau, \mathbf{u}}
\in X^{\an}$ lying on the generic point of $X$.

With the above notation, the assignment 
\begin{equation}
\label{eqn:intrinzic}
  \Theta_\tau=\sigma^{A}\times{\RR^B_{\geq 0}} \to \left\{ |\cdot|_{\tau, \mathbf{u}} \in X^{\an} \;\left|\; \mathbf{u} \in \Theta_\tau\right.\right\} \subset X^{\an}, \quad \mathbf{u} \mapsto  |\cdot|_{\tau, \mathbf{u}} 
\end{equation}
is a homeomorphism, cf. \cite[Proposition~3.4]{MN15}, \cite[4.2, 4.3]{GRW16}. We thus get an identification 
\begin{equation}
\label{eqn:coord:DS}
  \Theta_\tau \simeq \left\{ |\cdot|_{\tau, \mathbf{u}} \in X^{\an} \;\left|\; \mathbf{u} \in \sigma^{A} \times {\RR^{B}_{\geq 0}} \right.\right\} \subset X^{\an}.
\end{equation}
The polyhedron $\Theta_\tau$ is endowed with a natural integral affine structure.

Note that by the construction, 
$|\cdot|_{\tau, \mathbf{u}}$ is independent of the choice of  local equations $\z_v = 0$ and $\w_\rho = 0$ for $\Xscr_{\kappa, v}$ and $H_\rho$.  Thus, the embedding of $\Theta_\tau$ in the Berkovich space given in 
\eqref{eqn:intrinzic} is intrinsic, cf. \cite[p.~169]{GRW16}. 

Let $S_\sigma$ be another special stratum such that the Zariski closure $\overline{S_\sigma}$ of $S_\sigma$ contains $S_\tau$. 
Then,
$\Theta_{\sigma}$ is a face of $\Theta_\tau$. Moreover, the embedding of $\Theta_{\sigma}$ and $\Theta_\tau$ in $X^{\an}$ are compatible, as we now explain.
Since $\overline{S_\sigma} \supseteq S_\tau$ and $S_\sigma$ is a special stratum,
there exist a subset $A^\prime$ of 
$A$ and a subset $B^\prime$ of $B$ 
such that $S_\sigma$ is defined in $\Xscr$ at its generic point by 
$\z_{v} = 0$ for all $v \in A^\prime$ and 
$\w_{\rho} = 0$ for all $\rho \in B^\prime$. The extended simplex 
$\Theta_\sigma=\sigma^{A^\prime}\times {\RR^{B^\prime}_{\geq 0}}$
is the face of $\Theta_\tau=\sigma^{A}\times {\RR^{B}_{\geq 0}}$ consisting of those points with $u_v=0$ for $v\in A \setminus A^\prime$ and $u_\rho=0$ for $\rho\in B\setminus B^\prime$. 

Then, by the very definition of $|\cdot|_{\tau, \mathbf{u}}$, we see that 
$|\cdot|_{\tau, \mathbf{u}} = |\cdot|_{\sigma, \mathbf{u}}$ for any $\mathbf{u} \in \Theta_\sigma=\sigma^{A^\prime} \times \RR^{B^\prime}_{\geq 0}$. 
It follows that the embedding of 
$\Theta_{\sigma}$ in $X^{\an}$ coincides with the restriction to the face $\Theta_\sigma\subseteq \Theta_\tau$ of the embedding of $\Theta_\tau$ in $X^{\an}$.

\subsubsection*{Skeleton $\Sk(\Xscr, \cH)$}
The {\em skeleton} $\Sk(\Xscr, \cH) \subseteq X^{\an}$ associated with  
the strictly semistable pair $(\Xscr, \cH)$ is defined by  
\[
  \Sk(\Xscr, \cH) \colonequals \bigcup_{\tau \in \Delta(\Xscr, \cH)} \Theta_\tau.
\]
The skeleton can be identified with the support of  $\Delta(\Xscr, \cH)$, embedded in $X^{\an}$.

In the following, we often denote by $\sigma, \tau, \delta$ the elements of $\Delta(\Xscr, \cH)$. 
For two elements $\tau, \delta$ in $\Delta(\Xscr, \cH)$, we write 
\begin{equation}
\label{eqn:succeq}
\tau \succeq\delta
\end{equation}
if $\Theta_\delta$ is a face of $\Theta_\tau$, equivalently, if 
$\overline{S_\delta} \supseteq S_\tau$. 

\subsubsection*{Parallel rays in $\Delta(\Xscr, \cH)$ and parallel half-lines in $\Sk(\Xscr, \cH)$} We define parallel rays of $\Delta(\Xscr, \cH)$ as in Section~\ref{subsubsec:abstract:polyherdal:complex}. 
Identifying $\Sk(\Xscr, \cH)$ with the support of $\Delta(\Xscr, \cH)$, we define the notion of parallel half-lines for the skeleton.

\subsection{Tropicalization and the tropical function field of $(\Xscr, \cH)$}
We retain the notation and assumptions introduced in the previous subsection.

 We define $\Rat(\Delta(\Xscr, \cH))$ and $\ss \Rat(\supp{\Delta(\Xscr, \cH)})$ as in Section~\ref{subsubsec:abstract:polyherdal:complex}. The latter, $\ss \Rat(\supp{\Delta(\Xscr, \cH)})$, is referred to as the \emph{tropical function field of $(\Xscr, \cH)$}.

Let $f \in \Rat(\Xscr) = \Rat(X)$ be a nonzero element. Let $x \in \Sk(\Xscr, \cH)$, and let $\tau \in \Delta(\Xscr, \cH)$ and $\mathbf{u}\in \sigma^A \times \RR_{\geq0}^B$ such that $x =  |\cdot|_{\tau, \mathbf{u}}$, see \eqref{eqn:coord:DS}. 
We set 
\[
\trop(f)(x) \colonequals -\log|f|_{\tau, \mathbf{u}} \in \RR.
\]
This yields a function 
\[
\trop(f)\colon \Sk(\Xscr, \cH) \to \RR.
\]
We extend the definition by setting $\trop(0)=\infty$.
\begin{Proposition}
\label{prop:tropf:belongs:RXH}
For each $f \in \Rat(\Xscr) = \Rat(X)$, we have $\trop(f) \in \ss\Rat(\supp{\Delta(\Xscr, \cH)})$. 
\end{Proposition}

\begin{proof}
We may assume that $f \neq 0$. 
Take a special stratum $S_\tau$, $\tau\in \Delta(\Xscr, \cH)$. Let $A \subseteq V_\f$ and 
$B \subseteq \ss V_\infty$ be the subsets such that 
$S_\tau$ is a connected component of $Z_A \cap W_B$. We write $f = f_1/f_2$ with 
$f_1, f_2 \in \widehat{\OO}_{\Xscr, \xi_\tau}$, 
where $\xi_\tau$ is the generic point of $S_\tau$.

Using the notation from \eqref{def:mathu:S}, we represent $f_1$ and $f_2$ by admissible expansions   
\[
f_1= \sum_{\mathbf{m} \in \ZZ_{\geq 0}^{A\cup B}} a_{\mathbf{m}} \prod_{v\in A}\z_{v}^{m_v}  \prod_{\rho\in B}\w_{\rho}^{m_{\rho}}, 
\quad
f_2 = \sum_{\mathbf{n} \in \ZZ_{\geq 0}^{A \cup B}} a_{\mathbf{n}} \prod_{v\in A}\z_{v}^{n_v} \prod_{\rho\in B}\w_{\rho}^{n_{\rho}}. 
\]
For each $\mathbf{u}\in \Theta_\tau=\sigma^A\times \RR_{\geq 0}^B$, we have 
\begin{align*}
\trop(f)(\mathbf{u}) & = 
\min_{\mathbf{m}  \in  \ZZ_{\geq 0}^{A\cup B}, \; 
a_{\mathbf{m}}  \neq 0} \left\{ 
\sum_{v\in A} m_v u_{v}  + \sum_{\rho\in B} m_\rho u_\rho  \right\}
   - 
\min_{\mathbf{n} \in  \ZZ_{\geq 0}^{A \cup B}, \;
+ a_{\mathbf{n}}  \neq 0} \left\{ \sum_{v\in V} n_v u_{v} + \sum_{\rho\in B} n_\rho u_\rho\right\}. 
\end{align*}
It follows that $\trop(f)$ is a continuous piecewise integral $\QQ$-affine function. 
 
Furthermore, let $\varrho = \mathbf{u}+\RR_{\geq 0} \e_\rho$ be a ray in $\Theta_\tau$ parallel to $\rho \in B$. Let $W_\rho$ be the horizontal irreducible component of $\cH$ corresponding to $\rho$. The slope of $\trop(f)$ near infinity 
along $\varrho$ equals to the order of vanishing of $f$ along $W_\rho$, which remains constant for rays parallel to $\rho$. 

More generally, 
on each $\Theta_\tau$, $\trop(f)$ is the difference of the two functions $H_1, H_2$, where each $H_i$ is the tropicalization of $f_i$. Since $H_1$ is the minimum of affine functions consisting of a bounded term ($\sum_{v\in A} m_v u_{v}$) and an unbounded term ($\sum_{\rho\in B} m_\rho u_\rho$), and similarly for $H_2$,  
the asymptotic behavior near infinity of $H_i$ is captured by the minimum of their unbounded term 
\[
\ss H_{1, \infty} \colonequals \min_{\mathbf{m}  \in  \ZZ_{\geq 0}^{A\cup B}} \;\left\{\sum_{\rho\in B} m_\rho u_\rho\right\}, 
\quad 
\ss H_{2, \infty} \colonequals \min_{\mathbf{n}  \in  \ZZ_{\geq 0}^{A\cup B}} \;\left\{\sum_{\rho\in B} n_\rho u_\rho\right\}. 
\]

Let $S_\eta$ be the stratum of $\cY$ which contains $S_\tau$ and is a connected component of the intersection $\bigcap_{\rho\in B} W_\rho$. Denote by $\xi_\eta$ the generic point of $\eta$. 
The admissible expansion for $f_i$ at $\xi_\tau$ given above provides an admissible expansion of $f_i$ at $\xi_\eta$ of the form
\[
f_1= \sum_{\mathbf{m} \in \ZZ_{\geq 0}^{A\cup B}} \widehat{a}_{\mathbf{m}} \prod_{\rho\in B}\w_{\rho}^{m_{\rho}}, 
\quad
f_2 = \sum_{\mathbf{n} \in \ZZ_{\geq 0}^{A \cup B}} \widehat a_{\mathbf{n}} \prod_{\rho\in B}\w_{\rho}^{n_{\rho}},
\]
where $\widehat a_{\mathbf{m}} =  a_{\mathbf{m}} \prod_{v\in A}\z_{v}^{m_v} $ and $\widehat a_{\mathbf{n}}= a_{\mathbf{n}} \prod_{v\in A}\z_{v}^{n_v}$.
The associated cone $\Theta_\eta$ of $\Upsilon$ lives in the Berkovich analytification of $X$, when $K$ is viewed as a trivially valued field. The tropicalization of $f_i$ on $\Theta_\eta$ is given by the admissible expansions in the form $\trop(f_i) = \min_{\mathbf{m}  \in  \ZZ_{\geq 0}^{A\cup B}} \;\left\{\sum_{\rho\in B} m_\rho u_\rho\right\}$ and $\min_{\mathbf{n}  \in  \ZZ_{\geq 0}^{A\cup B}} \;\left\{\sum_{\rho\in B} n_\rho u_\rho\right\}$, and does not depend on the choice of admissible expansions (see \cite[Proposition~3.1]{JM12}). 
This is equal to $H_{i, \infty}$.

By the definition of parallel half-lines in $\supp{\Delta(\Xscr, \cH)}$, we infer that $\trop(f)$ has the same slope near infinity along two 
parallel half-lines in $\supp{\Delta(\Xscr, \cH)}$. 

We infer that $\trop(f)$ belongs to $\ss\Rat(\supp{\Delta(\Xscr, \cH)})$. 
\end{proof}

Let $\bar{K}$ be an algebraic closure of $K$. Let 
$X_{\bar{K}} = X \times_K  \bar{K}$. Recall that 
we assume that $X$ and each connected component of the intersection of any collection of irreducible components of $\cH$ are geometrically connected. 

\begin{Proposition}
\label{prop:tropf:belongs:RXH:2}
The tropicalization map extends to $\trop\colon \Rat(X_{\bar{K}})$ $\to \ss\Rat(\supp{\Delta(\Xscr, \cH)})$. 
\end{Proposition}

\begin{proof} Let $K^\prime$ be a finite field extension of $K$, and let $f \in \Rat(X_{K^\prime})$, 
where $X_{K^\prime} \colonequals X \times_K K^\prime$. 
Denote by $R^\prime$ the integral closure of $R$ in $K^\prime$. Then, $R^\prime$ is a discrete valuation ring (see \cite[Chap.~12, Th.~72]{M80}). 
Let $\varpi^\prime$ be a uniformizer of $R^\prime$, and recall that $\varpi$ denotes a uniformizer of $R$. 
We can write 
\[
\varpi = \varpi^{\prime e} u^\prime \qquad \text{with } e \geq 1 \text{ and } u^\prime \in R^{\prime\times}.
\]
The base change $\Xscr \times_{R} R^\prime$ is expressed \'etale locally as 
\[
\Spec \bigl(R^\prime[\x_0, \ldots, \x_d]/(\x_0 \cdots \x_r - \varpi^{\prime e} u^\prime)\bigr).
\]

By Theorem~\ref{thm:Cartwright}, there exists a regular strictly semistable model $\Xscr^\prime \to \Spec(R^\prime)$ along with a surjective morphism $\mu\colon \Xscr^\prime \to \Xscr \times_{R} R^\prime$. Let $\pi^\prime\colon \Xscr \times_{R} R^\prime \to \Xscr$ be the natural morphism and set $\pi \colonequals \pi^\prime \circ \mu 
\colon \Xscr^\prime \to \Xscr$. 
Let $\cH^\prime$ be the strict transform of $\cH$ by $\pi$. Then, $(\Xscr^\prime, \cH^\prime)$ is a strictly semistable pair on $X_{K^\prime}$, and  $\Delta(\Xscr^\prime, \cH^\prime)$ forms an integral $\QQ$-affine refinement of $\Delta(\Xscr, \cH)$. In particular, the skeleton $\Sk(\Xscr, \cH)$ is naturally identified with $\supp{\Delta(\Xscr^\prime, \cH^\prime)}$. By our assumption, the recession fan of $\Delta(\Xscr^\prime, \cH^\prime)$ is the same as that of $\Delta(\Xscr, \cH)$. 

By applying Proposition~\ref{prop:tropf:belongs:RXH} to $(\Xscr^\prime, \cH^\prime)$, we obtain the tropicalization map 
\[
\trop\colon \Rat(X_{K^\prime}) \to \ss\Rat(\supp{\Delta(\Xscr^\prime, \cH^\prime)}) = \ss\Rat(\supp{\Delta(\Xscr, \cH)}). 
\]
For $f \in \Rat(X)$, viewed in $\Rat(X_{K^\prime})$, the tropicalization $\trop(f)$ obtained above coincides with the previously defined tropicalization. 

Hence, we obtain the tropicalization map 
$\trop\colon \Rat(X_{\bar{K}})$ $\to \ss\Rat(\supp{\Delta(\Xscr, \cH)})$. 
\end{proof}

\subsection{Faithful tropicalization}
We recall the definition of faithful tropicalization; see~\cite[\S9]{GRW16} for further details. 

Let $(\Xscr, \cH)$ be a strictly semistable pair, where $\Xscr$ is a strictly semistable model of $X$ defined over $K$. 

Given nonzero rational functions $f_1, \dots, f_n \in \Rat(X_{\bar K})$, let $U$ be the open subset of $X_{\bar K}$ obtained by removing the support of the divisors $\mathrm{div}(f_j)$, $j=1, \dots, n$. Consider the map $\phi=(f_1, \dots, f_n) \colon U \to {\mathbf{\ss G_{m,\bar K}^n}}$ from $U$ to the algebraic torus ${\mathbf{\ss G_{m,\bar K}^n}}$. The tropicalization of $\phi$ induces a map $\ss\phi^{\trop} \colon U^{\an} \to \RR^n$. 

Since the skeleton $\Sk(\Xscr,\cH)$ is included in $U^{\an}$, the restriction of $\ss\phi^\trop$ gives a map 
\[
\ss\phi^{\trop} \colon \Sk(\Xscr,\cH) \to \RR^n.
\]
We say that $\ss\phi^\trop$ is \emph{faithful} if the following conditions are satisfied:
\begin{itemize}
    \item The restriction of $\ss\phi^\trop$ to the skeleton is injective. 
    \item On each $\Theta_\tau$, $\tau \in \Delta(\Xscr, \cH)$, the restriction $\rest{\ss\phi^\trop}{\Theta_\tau}$ is piecewise integral $\QQ$-affine.
    \item The inverse of $\rest{\ss\phi^\trop}{\Theta_\tau}$ from $\ss\phi^\trop(\Theta_\tau)$ to $\Theta_\tau$ is also piecewise integral $\QQ$-affine. 
\end{itemize}
Equivalently, $\ss\phi^\trop$ is faithful if it is  injective, its image is an integral $\QQ$-affine polyhedral complex in $\RR^n$, and the inverse map from the image of $\ss\phi^\trop$ to $\Sk(\Xscr,\cH)$ is piecewise integral $\QQ$-affine.  

\subsection{Semifields and finite generation}
\label{subsec:tropical:semifield}
We set 
\[
\TT \colonequals 
\RR \cup \{\infty\}. 
\]
For $a, b \in \TT$, we set 
\[
a \oplus b \colonequals \min\{a, b\}, \quad 
a \odot b \colonequals a + b. 
\]
With these operations, $\TT$ forms a semifield, where $\infty$ is the identity element for $\oplus$ and $0$ is the identity element for $\odot$. We call $\TT$ the {\em tropical semifield}.  

Similarly, we define the tropical subsemifield $\TT\QQ$ as 
\[
\TT\QQ \colonequals \QQ \cup \{\infty\} \subseteq \TT. 
\]

Let $S$ be a semiring with binary operations $\oplus$ and $\odot$. We say that $S$ is a {\em semiring  over $\TT\QQ$} if there exists an injective semiring homomorphism 
$\TT\QQ \to S$. In this case, 
for any $a \in \TT\QQ$ and $s \in S$, the operation $a \odot s$ is well-defined via this homomorphism. Moreover, if $S$ is a semifield, we call it a semifield {\em over $\TT\QQ$}. 

\begin{Example}[$\ss\Rat(\supp\Sigma)$]
Let $\ss\Rat(\supp\Sigma)$ be defined as in Definition~\ref{def:rat:Delta}. For $F, G \in \ss\Rat(\supp\Sigma)$ and $a \in \TT\QQ$, we define
\[
F \oplus G \colonequals 
\min\{F, G\}, \quad 
F \odot G \colonequals 
F + G, \quad 
a \odot F \colonequals 
a + F. 
\]
Under these operations,  $F \oplus G$, $F \odot G$, and $a \odot F$ all belong to $\ss\Rat(\supp\Sigma)$. Moreover,
$\ss\Rat(\supp\Sigma)$ forms a semiring over $\TT\QQ$.  Furthermore, for any $F \in \ss\Rat(\supp\Sigma)$ other than $\infty$, we have 
$F \odot (-F) = 0$. This shows that $\ss\Rat(\supp\Sigma)$ is a semifield over $\TT\QQ$.  
\end{Example}

Let $S$ be a semiring over $\TT\QQ$, and $s_1, \ldots, s_r \in S$. Define the subsemiring of $S$ generated by $s_1, \dots, s_r$ as
\[
\TT\QQ[s_1, \ldots, s_r]
\colonequals 
\left\{\left. \bigoplus_{(i_1, \ldots, i_r) \in A \subseteq \ZZ_{\geq 0}^r} a_{i_1,\ldots, i_r} \odot s_1^{\odot i_1}\odot\cdots\odot s_r^{\odot i_r} \;\right|\; A \textrm{ finite and } a_{i_1,\ldots, i_r} \in \TT\QQ\right\}
\]
The semifield of fractions of $\TT\QQ[s_1, \ldots, s_r]$ denoted $\TT\QQ(s_1, \ldots, s_r)$ is defined by localizing at its nonzero elements. If $S$ is a semifield, then $\TT\QQ(s_1, \ldots, s_r) \subseteq S$.

\begin{Definition}[Finitely generated semifields]
\label{def:fin:gen}
Let $S$ be a semifield over $\TT\QQ$. We say that $S$ is {\em finitely generated as a semifield over $\TT\QQ$} if 
there exist $s_1, \ldots, s_r \in S$ such that 
$S = \TT\QQ(s_1, \ldots, s_r)$. 
\end{Definition}

\subsection{Congruences}
\label{subsec:congruences}
Let $S$ be a semiring with binary operations $+$ for addition and $\cdot$ for multiplication. An equivalence relation $\sim$ on $S$ is a {\em congruence} if it satisfies the following property: 
For any $a, b, a^\prime, b^\prime \in S$, if $a \sim b$ and $a^\prime \sim b^\prime$, then 
\[
a + a^\prime \sim b + b^\prime \qquad \text{and} \qquad a \cdot a^\prime \sim b \cdot b^\prime.
\]
We refer to~\cite[\S I.7]{HW93} for further details.

If $\sim$ is a congruence on $S$, then the quotient $\overline{S} \colonequals S/\!\sim$ is itself a semiring. Conversely, if $\varphi\colon S \to S^\prime$ is a semiring homomorphism, then the kernel 
\[
{\rm Ker}(\varphi) \colonequals 
\left\{(s_1, s_2) \in S^2 \mid \varphi(s_1) = \varphi(s_2)\right\}
\]
is a congruence on $S$. In the context of semirings, the operation of taking the quotient of a ring by an ideal is replaced by taking the quotient by a congruence. 

\begin{Definition}[Finite generation property for congruences]
Let $\sim$ be a congruence on $S$. We say that 
$\sim$ is {\em finitely generated} if 
there exist elements $F_1, \ldots, F_r, G_1, \ldots, G_r \in S$ with 
$F_i \sim G_i$ for all $i$ such that the following holds:

If $\sim^\prime$ denotes the minimal congruence on $S$ generated by the relations $F_1 \sim^\prime G_1, \ldots, F_r \sim^\prime G_r$, then $\sim^\prime$ coincides with $\sim$. 
\end{Definition}


\section{Characterization theorem}
\label{sec:proof:Q1}
As before, let $K$ be a complete discrete valuation field. Given a finite extension $K^\prime$ of $K$,  we denote by $R^\prime$ the integral closure of $R$ in $K^\prime$, and by $\kappa^\prime$ the residue field of $R^\prime$. 

Let $X$ be a smooth, geometrically connected, projective variety over $K$, and suppose that $X$ admits a strictly semistable model $\Xscr$ over $R$.
Let $(\Xscr, \cH)$ be a strictly semistable pair, and assume that each connected component of the intersection of any collection of irreducible components of $\cH$ is geometrically connected. Let $\Delta(\Xscr, \cH)$ be the dual complex of $(\Xscr, \cH)$, and let $\Sk(\Xscr, \cH) \subseteq X^{\an}$ be the associated skeleton. Denote by $\Upsilon$ the recession fan of $\Delta(\Xscr, \cH)$.

In this section, we prove Theorem~\ref{thm:Q:1}. Specifically, we show that for any $\Phi$ in the tropical function field $\Rat(\supp{\Delta(\Xscr, \cH)})$, there exists an element $\phi \in \Rat(X_{\bar{K}})$ such that $\trop(\phi) = \Phi$.

To establish this assertion, we need the following result, 
see~\cite{Hart01,Cartwright}.

\begin{Theorem}
\label{thm:Cartwright}
Notation as above:
\begin{enumerate}[label=\textup{(C\arabic*)}]
\item \label{Car1}
There exist a finite extension field $K^{\prime}/K$, a strictly semistable model $\Xscr^\prime \to \Spec(R^\prime)$ of $X^\prime \colonequals X \times_{K} K^\prime$ and a morphism $\pi\colon \Xscr^\prime \to \Xscr$ with the following properties\textup{:} 
\begin{enumerate}
\item[\textup{(i)}]
Let $\cH^\prime$ be the strict transform of $\cH$ by $\pi$. Then, $(\Xscr^\prime, \cH^\prime)$ is a strictly semistable pair, and 
$\Delta(\Xscr^\prime, \cH^\prime)$ is a subdivision of $\Delta(\Xscr, \cH)$ with the same recession fan. 
\item[\textup{(ii)}]
Let $Z_v^\prime$, $v\in V_\f^\prime$, be the irreducible components of $\Xscr^\prime_{\kappa^\prime}$, and $W_\rho^\prime$, $\rho\in \ss V_\infty$ be the irreducible components of $\cH^\prime$. Then, for any subset $\emptyset\neq A\subseteq V_\f^\prime$ and $B\subseteq \ss V_\infty$, the intersection $\bigcap_{v\in A} Z_v^\prime \cap \bigcap_{\rho\in B} W_\rho^\prime$ is either empty or irreducible. 
\end{enumerate}
\item \label{Car2}
Let $\Phi \in \Rat(\Sk(\Xscr, \cH))$. 
Then, there exists a finite extension field $K^{\prime}/K$ and a strictly semistable model  $\Xscr^\prime \to \Spec(R^\prime)$ of $X^\prime \colonequals X \times_{K} K^\prime$ satisfying Properties~\textup{(i)} and \textup{(ii)} of \ref{Car1} such that 
$\Phi \in \Rat(\Delta(\Xscr^\prime, \cH^\prime))$. 
\end{enumerate}
\end{Theorem}

\begin{proof} 
We follow the desingularization procedure described in~\cite{Hart01}, which involves performing successive below-ups, adapted to the strata obtained from intersections of the components of $\cH$ and those of the special fiber. This process results in a semistable model $\Xscr^\prime \to \Spec(R^\prime)$ with a map $\pi\colon \Xscr^\prime \to \Xscr$. The corresponding dual complex yields a triangulation of the dual complex of $(\Xscr, \cH)$. More precisely, it forms an order $m$ subdivision of $\Delta(\Xscr, \cH)$, in the sense of~\cite{Cartwright}, where $m$ is the ramification index of $K^\prime/K$. For $m>1$, this subdivision has no multiple faces. Moreover, by our assumption on $\cH$, the recession fan does not change, from which \ref{Car1} follows.

For each $\Phi \in \Rat(\Sk(\Xscr, \cH))$, there exists a large $m$ and an $m$-subdivision of $\Delta(\Xscr, \cH)$ as above over which $\Phi$ is facewise integral $\QQ$-affine. The second statement follows.
\end{proof}

Our strategy for proving Theorem~\ref{thm:Q:1} goes as follows. 
By replacing $K$ with $K^\prime$ and $\Xscr$ with $\Xscr^\prime$, we may and do assume that our model $\Xscr$ satisfies the properties stated in~\ref{Car1} and~\ref{Car2} of Theorem~\ref{thm:Cartwright}, and that 
$\Phi$ is facewise integral $\QQ$-affine on each (extended simplex) $\Theta_\tau$ associated with a special stratum $S_\tau$ of $\Xscr$. Moreover, 
passing to a large field extension $K^\prime$ of $K$ if necessary, we can assume that the constant terms determined by the restriction of $\Phi$ to $\Theta_\tau$ (see proof of Proposition~\ref{prop:construction:g:sigma}) belong to the value group of $K$, for all $\tau$, and that the residue field $\kappa$ has sufficiently large cardinality (if it is finite).

Recall that $\Sk(\Xscr, \cH) = \bigcup_{\tau \in \Delta(\Xscr, \cH)} \Theta_\tau$ (see Section~\ref{subsec:semistable:models}). Given a function $F \in \ss\Rat(\Sk(\Xscr, \cH))$ and a pair $\tau \preceq \eta$ in $\Delta(\Xscr, \cH)$, we denote by $\rest F\tau$ and $\rest F{\eta\setminus\tau}$ the restrictions of $F$ to $\Theta_\tau$ and $\Theta_\eta \setminus \Theta_\tau$, respectively.

For each $\tau \in \Delta(\Xscr, \cH)$, we will construct two functions with {\em opposite} behaviors. The first function denoted by $f_\tau \in \Rat(\Xscr) = \Rat(X)$ satisfies 
\begin{enumerate}
\item[(i)]
$\rest{\trop(f_\tau)}{\tau} \equiv 0$; 
\item[(ii)]
$\rest{\trop(f_\tau)}{\eta \setminus \tau}  > 0$ for any $\eta \in \Delta(\Xscr, \cH)$. Furthermore, if $\eta \in \Delta(\Xscr, \cH)$,   
then the slope of 
$\rest{\trop(f_\tau)}{\eta\setminus\tau}$ near infinity along 
any half-line of $\eta\setminus \tau$ that is not parallel to $\tau$ is positive.
\end{enumerate}

We will show from the construction that $\rest{\trop(f_\tau)}{\eta\setminus\tau}$ is concave and the slope of $\rest{\trop(f_\tau)}{\eta\setminus\tau}$ near infinity along any ray of $\eta\setminus \tau$ that is not parallel to a ray of $\tau$ is positive, from which we deduce the more general statement about half-lines. 

\smallskip
The second function denoted by 
$g_\tau \in \Rat(\Xscr) = \Rat(X)$ satisfies 
\[
\rest{\trop(g_\tau)}{\tau} = \rest{\Phi}{\tau}. 
\]

With these two functions $f_\tau, g_\tau \in 
\Rat(\Xscr) = \Rat(X)$, 
we set  
\begin{equation}
\label{eqn:construction:h}
\phi\colonequals 
\sum_{\tau \in \Delta(\Xscr, \cH)}
\beta_\tau f_\tau^{N_\tau} g_\tau \in \Rat(\Xscr) = \Rat(X), 
\end{equation}
where the coefficients $\beta_\tau \in R^\times$  and the integers $N_\tau \gg 1$ will be suitably chosen. We then show that $\trop(\phi) = \Phi$, as desired.

\subsection{Construction of $f_\tau$}
By the definition of a strictly semistable pair, 
$\Xscr$ is covered by open subsets that admit an \'etale morphism $\mathscr{U} \to \Spec(R[\x_0, \ldots, \x_d]/\langle \x_0 \cdots \x_r - \varpi\rangle)$. 

Let $\tau$ be an element of $\Delta(\Xscr, \cH)$, and let $S_\tau$ be the special stratum of $(\Xscr, \cH)$ corresponding to $\tau$ defined by $A_\tau \subseteq V_\f$ and $B_\tau \subseteq \ss V_\infty$. Let $\xi_\tau$ denote the generic point of $S$.

By \ref{Car1} in Theorem~\ref{thm:Cartwright}, we have 
\[
S_\tau = \bigcap_{v \in A_\tau} Z_v \cap \bigcap_{\rho \in B_\tau} W_\rho.
\]

Recall that $\cY = \Xscr_\kappa + \cH$. To simplify the notation, we set
\[
D_{\tau} \colonequals \sum_{v \in A_\tau} Z_v + \sum_{\rho \in B_\tau} W_\rho, \qquad 
E_{\tau} \colonequals 
\cY - D_\tau = 
\sum_{v \in V_\f \setminus A_\tau} Z_{v} + \sum_{\rho \in \ss V_\infty \setminus B_\tau} W_\rho.  
\]

We pick an ample line bundle $\Lscr$ over $\Xscr$.

\begin{Lemma}
\label{lemma:fsigma:1}
For a sufficiently large positive integer $N$, there exists a section 
\[
s \in H^0\Bigl(D_{\tau}, \rest{\Lscr^{\otimes N}\bigl(-E_{\tau}\bigr)}{D_{\tau}}\Bigr)
\] 
such that $s(\xi_\tau) \neq 0$. 
\end{Lemma}

\begin{proof}
Since $\Lscr$ is ample, the restriction $\rest{\Lscr}{D_{\tau}}$ is also ample. 
Therefore, for sufficiently large $N$, $\rest{\Lscr^{\otimes N}\bigl(-E_{\tau}\bigr)}{D_{\tau}}$ is 
globally generated. It follows that there exists $s \in  H^0\Bigl(D_{\tau}, \rest{\Lscr^{\otimes N}\bigl(-E_{\tau}\bigr)}{D_{\tau}}\Bigr)$ such that 
$s(\xi_\tau) \neq 0$. 
\end{proof}

\begin{Lemma}
\label{lemma:fsigma:2}
Let $s$ be as in Lemma~\ref{lemma:fsigma:1}. Then, there exists a section
\[
s^\prime \in 
H^0\left(\cY, \rest{\Lscr^{\otimes N}}{\cY}\right)
\]
such that $\rest{s^\prime}{D_\tau} = s$ and $\rest{s^\prime}{E_\tau} = 0$. 
\end{Lemma}

\begin{proof}
The restriction map $\OO_{\cY} \twoheadrightarrow \OO_{E_\tau}$ induces, after tensoring by $\Lscr^{\otimes N}$, the exact sequence  
\[
0 \to \rest{\Lscr^{\otimes N}(-E_\tau)}{D_\tau} \to \rest{\Lscr^{\otimes N}}{\cY} \to 
\rest{\Lscr^{\otimes N}}{E_\tau} \to 0. 
\] 
In the corresponding long exact sequence of cohomologies, let $s^\prime \in H^0\left(\cY, \rest{\Lscr^{\otimes N}}{\cY}\right)$ be the push-out of $s \in H^0\left(D_\tau, \rest{\Lscr^{\otimes N}(-E_\tau)}{D_\tau}\right)$. 
In other words, $s^\prime$ is the zero extension of $s$ to $\cY$. 
Then, $s^\prime$ has the desired properties.
\end{proof}

\begin{Lemma}
\label{lemma:fsigma:3}
If $N$ is a sufficiently large positive integer, then the restriction map 
\[
H^0(\Xscr, \Lscr^{\otimes N}) \to H^0(\cY, \rest{\Lscr^{\otimes N}}{\cY})
\]
is surjective. For such an $N\gg 1$, 
we take $s$ as in Lemma~\ref{lemma:fsigma:1} and  
$s^\prime$ as in Lemma~\ref{lemma:fsigma:2}. Then, there exists 
\[
\widetilde{s} \in 
H^0\left(\Xscr, \Lscr^{\otimes N}\right)
\]
such that $\rest{\widetilde{s}\,}{\cY} = s^\prime$. 
\end{Lemma}

\begin{proof}
The first assertion follows from the vanishing $H^1\left(\Xscr, \Lscr^{\otimes N}(-\cY)\right) = 0$ for sufficiently large $N$. The second statement follows directly from the first. 
\end{proof}

\begin{Proposition}
\label{prop:construction:f:sigma}
Notation as above, for each $\tau \in \Delta(\Xscr, \cH)$, there exists $f_{\tau} \in \Rat(\Xscr) = \Rat(X)$ such that 
\[
\rest{\trop(f_{\tau})}{\tau} \equiv 0\quad \text{and} \quad \rest{\trop(f_{\tau})}{\eta\setminus \tau} >0\quad \text{for each}\,\, \eta \in \Delta(\Xscr, \cH).
\]
Furthermore, $\rest{\trop(f_{\tau})}{\eta}$ is concave for any 
$\eta \in \Delta(\Xscr, \cH)$, and 
the slope of $\rest{\trop(f_{\tau})}{\eta\setminus \tau}$ near infinity along any ray of $\eta$ that is not parallel to a ray of $\tau$ is positive. 
\end{Proposition}

\begin{proof}
Let $\widetilde{s} \in H^0\left(\Xscr, \Lscr^{\otimes N}\right)$ 
be as in Lemma~\ref{lemma:fsigma:3}. 
For each $\eta \in \Delta(\Xscr, \cH)$, let $\xi_\eta$ denote the generic point 
of stratum $S_\eta$ in $\cY$. Since $\Lscr$ is ample, for sufficiently large $N$, there exists $\widetilde{s}_0 \in H^0\left(\Xscr, \Lscr^{\otimes N}\right)$ such that 
$\widetilde{s}_0(\xi_\eta) \neq 0$ for all $\eta \in \Delta(\Xscr, \cH)$. Taking $N$ larger if necessary, we may and do assume that the integer $N$ is the same for both $\widetilde{s}$ and $\widetilde{s}_0$. 

We define 
\[
f_\tau \colonequals \widetilde{s}/\widetilde{s}_0 \in \Rat(\Xscr) = \Rat(X)
\]
and we aim to show that $f_\tau$ has desired properties. 

\smallskip

First, we verify that $\rest{\trop(f_\tau)}{\tau} \equiv 0$. Recall that $A_\tau \subseteq V_\f$ and $B_\tau \subseteq \ss V_\infty$ are the subsets with $S_\tau = \bigcap_{v \in A_{\tau}} Z_v \cap \bigcap_{\rho\in B_\tau} W_\rho$. We choose local equations $\z_v = 0$ for $Z_v$ and $\w_\rho=0$ for $W_\rho$ in $\Xscr$ at the generic point $\xi_\tau$ of $S_\tau$. 

Since $f_\tau(\xi_\tau) = \frac{\widetilde{s}(\xi_\tau)}{\widetilde{s}_0(\xi_\tau)} \neq 0$ by Lemma~\ref{lemma:fsigma:1}, for an admissible expansion of $f_\tau$ in $\widehat{\OO}_{\Xscr, \xi_\tau}$ 
\[
f_\tau = \sum_{\mathbf{m}=(m_v)_v\times (m_\rho)_\rho\in {\ZZ^{A_\tau\cup B_\tau}_{\geq 0}}} a_{\mathbf{m}}  \prod_{v\in A_\tau}\z_{v}^{m_v}\prod_{\rho\in B_\tau} \w_\rho^{m_{\rho}},
\]
with each coefficient $a_{\mathbf{m}}$ is either zero or a unit in $\widehat{\OO}_{\Xscr, \xi_\tau}$, we have 
$a_{\mathbf{0}} \neq 0$, and so we get 
\[
|f_\tau|_{\tau,\mathbf{u}} = 1\quad \text{for any} \,\,
\mathbf{u} \in \sigma^{A_\tau} \times {\RR^{B_\tau}_{\geq 0}}.
\]
 It follows that $\rest{\trop(f_\tau)}{\tau} \equiv 0$, proving the first claim. 

\smallskip

Next, we verify that $\rest{\trop(f_\tau)}{\eta\setminus\tau} > 0$ for each element $\eta$ of $\Delta(\Xscr, \cH)$. We may assume that $\eta \setminus\tau \neq \emptyset$, otherwise, there is nothing to prove. This implies that $\eta \not \preceq \tau$.

For the special stratum $S_\eta$ of $\cY$ corresponding to $\eta$ with generic point $\xi_\eta$, consider the corresponding subsets $A_{\eta} \subseteq V_\f$ and $B_\eta \subseteq \ss V_\infty$ with 
\[
S_\eta = \bigcap_{v\in A_\eta} Z_v \cap \bigcap_{\rho\in B_\eta} W_\rho.
\]
Let $\delta = \eta \cap \tau$, with $A_{\delta} = A_\eta \cap A_\tau$ and $B_{\delta} = B_\eta \cap B_\tau$. Since $\eta \not\preceq \tau$, we have $|A_{\delta}|+|B_{\delta}| < |A_{\eta}|+|B_{\eta}|$. 

We choose local equations $\z_v=0$ and $\w_\rho = 0$ at $\xi_\eta$ in $\Xscr$  for $Z_v$ and $W_\rho$, respectively, for $v\in A_\eta, \rho\in B_\eta$.

Let $v\in A_\eta \setminus A_\tau$. 
Since $v \not\in A_\tau$, $Z_v$ does not appear in $D_{\tau}$. By the construction of $\widetilde{s}$, we see that $\rest{\widetilde{s}}{Z_v} \equiv 0$. 
It follows that $\rest{\widetilde{s}}{(\z_v = 0)} \equiv 0$ for 
each $v \in A_\eta \setminus A_\tau$. 

Let now $\rho \in B_\eta \setminus B_\tau$. Since $\rho \not\in B_{\tau}$, $W_\rho$ does not appear in $D_{\tau}$. By the construction of $\widetilde{s}$, we see that $\rest{\widetilde{s}}{W_\rho} \equiv 0$.  It follows that $\rest{\widetilde{s}}{(\w_\rho = 0)} \equiv 0$ for 
each $\rho \not\in B_\tau$

Identifying $\Theta_{\eta} = \sigma^{A_\eta} \times \RR^{B_\eta}_{\geq 0}$ and  $\Theta_{\delta} = \sigma^{A_\delta} \times \RR^{B_\delta}_{\geq 0}$, the inclusion $\Theta_\delta \subseteq \Theta_\eta$ is given by 
\[
\Theta_\delta = \Bigl\{\mathbf{u} = (u_v)_{v\in A_\eta}\times (u_\rho)_{\rho \in B_\eta} \in \Theta_\eta \,\bigl |\, u_v=u_\rho = 0 \text{ for all } v\in A_\eta \setminus A_\delta, \rho \in B_\eta \setminus B_\delta\Bigr\}.
\]

Since $\widetilde{s}_0(\xi_\eta) \neq 0$ 
and $\rest{\widetilde{s}}{(\z_{v} = 0)}$ and $\rest{\widetilde{s}}{(\w_{\rho} = 0)}$ both vanish for $v\in A_\eta \setminus A_\delta$ and $\rho\in B_\eta \setminus B_\delta$, the admissible expansion of $f_\tau$ in $\widehat{\OO}_{\Xscr, \xi_\eta}$ 
\begin{equation}
\label{eqn:adm:ftau:eta}
f_\tau = \sum_{\mathbf{m} = (m_v)_v\times (m_\rho)_{\rho} \in {\ZZ^{A_\eta \cup B_\eta}_{\geq 0}}} a_{\mathbf{m}} \prod_{v\in A_\eta}\z_{v}^{m_v}\prod_{\rho\in B_\eta} \w_{\rho}^{m_{\rho}} 
\end{equation}
verifies $m_{v} > 0$ and $m_\rho>0$ for all $v\in A_\eta \setminus A_\delta$ and $\rho\in B_\eta \setminus B_\delta$.
It follows that 
$|f_\tau|_{\eta,\mathbf{u}} < 1$ if $\mathbf{u} \in \Theta_\eta \setminus \Theta_\delta$. We conclude that  $\rest{\trop(f_\tau)}{\eta\setminus\tau} > 0$, giving the second claim. 

It remains to show the last assertion. It follows from 
\eqref{eqn:adm:ftau:eta} that $\rest{\trop(f_\tau)}{\eta}$ is concave, being the minimum of a collection of affine functions.  

The rays of $\eta\setminus \tau$ that are not parallel to rays of $\tau$ correspond precisely to the elements of $B_\eta \setminus B_\tau$. Let $\rho$ be such an element. Since $\rest{\widetilde{s}}{(\w_\rho = 0)} \equiv 0$, we see that $m_\rho >0$ for all $\mathbf{m}$ with $a_{\mathbf m} \neq 0$. Therefore,  
the slope of $\rest{\trop(f_\tau)}{\eta}$ 
is positive along $\rho$, as required. 
\end{proof}

\begin{Corollary}
For any $\tau$ and $\eta$ in $\Delta(\Xscr, \cH)$, the slope of 
$\rest{\trop(f_\tau)}{\eta}$ near infinity along any half-line of $\eta$ that is not parallel to $\tau$ is positive. 
\end{Corollary}

\begin{proof}
By Proposition~\ref{prop:tropf:belongs:RXH}, $\trop(f_\tau) \in \Rat(\supp{\Delta(\Xscr, \cH)})$. 
The statement now follows from the concavity of $\rest{\trop(f_\tau)}{\eta}$ and positivity of the slope of $\rest{\trop(f_\tau)}{\eta}$ near infinity along any ray of $\eta$ that is not parallel to a ray of $\tau$. 
\end{proof}

\subsection{Construction of $g_\tau$}
In this subsection, we construct $g_\tau$. 
We need the following lemma. 

\begin{Lemma}
\label{lem:good:choice:local:eqns}
Let $E \subseteq \Xscr$ be an irreducible divisor. Given
points $x_1, \ldots, x_n$ on $E$ and points $y_1, \dots, y_m$ in $\Xscr \setminus E$, 
there is a rational function $u$ which vanishes on $E$ with order of vanishing one, belongs to each local ring $\OO_{\Xscr,x_i}$, $i = 1, \ldots, n$, and is a unit in $\OO_{\Xscr,y_j}$ for each $j=1, \dots, m$.
\end{Lemma}

\begin{proof}
We take an ample line bundle $\Lscr$ on $\Xscr$. 
There is a global section $s \in H^0(\Xscr, \Lscr^{\otimes N})$ for sufficiently large $N$ such that $s(x_i) \neq 0$ for any $i = 1, \ldots, n$ and $s(y_j) \neq 0$ for any $j = 1, \ldots, m$. Let $H = \zero(s)$. Then, $\Lscr \cong \OO_{\Xscr}(H)$, and $H$ does not pass through $x_1, \ldots, x_n$, $y_1, \ldots, y_m$ and does not contain $E$. We choose another large integer $N^\prime$ such that $N^\prime H - E$ is very ample. Taking a generic section $u \in H^0(\Xscr, \OO_{\Xscr}(N^\prime H - E))$, the corresponding rational function satisfies all the required properties. 
\end{proof}

We will apply this lemma to $E = Z_v$ for $v \in V_\f$ (resp. $E=W_\rho$ for $\rho\in \ss V_\infty$). The points $x_i$ will be the generic points $\xi_\tau$ of the special strata $S_\tau$ with $v\in \zeta(\tau)$ (resp. $\rho\in \zeta(\tau)$), while the points $y_j$ will be the generic points of the special strata $S_\eta$ with $v\not\in \zeta(\eta)$ (resp. $\rho\not\in \zeta(\eta)$).

\begin{Proposition}
\label{prop:construction:g:sigma}
Let $\Phi \in \Rat(\Delta(\Xscr, \cH))$.  Then, for each $\tau \in \Delta(\Xscr, \cH)$, 
there exists $g_\tau \in \Rat(\Xscr) = \Rat(X)$ such that 
$\trop(g_\tau)$ belongs to $\Rat(\Delta(\Xscr, \cH))$ 
and $\rest{\trop(g_\tau)}{\tau} = \rest{\Phi}{\tau}$. 
\end{Proposition}

Note that, in particular, $\trop(g_\tau)$ is integral $\QQ$-affine 
on each cell of $\Delta(\Xscr, \cH)$. 

\begin{proof}
  As before, let $S_\tau$ be the special stratum corresponding 
to $\tau$, and $A_\tau \subseteq V_\f$ and $B_\tau \subseteq \ss V_\infty$ be such that 
$S_\tau = \bigcap_{v \in A_\tau}Z_v \cap   \bigcap_{\rho \in B_\tau}W_\rho$. Applying Lemma~\ref{lem:good:choice:local:eqns}, we choose local equations  $\z_v = 0$ for $Z_v$ and $\w_\rho = 0$ for $W_\rho$ in $\Xscr$ at the generic point $\xi_\tau$ of $S_\tau$. with the property that $\z_v$ and $\w_\rho$ belong to all the local rings $\OO_{\Xscr, \xi_\eta}$ for $\eta \in \Delta(\Xscr, \cH)$, and they are unit if $v$ is not a vertex and $\rho$ is not a ray of $\eta$.

Identifying $\Theta_\tau$ with $\sigma^{A_\tau}\times \RR^{B_\tau}_{\geq 0}$, the restriction $\rest{\Phi}\tau$ can be written in the form 
\begin{equation}\label{eqn:pullback:Phi}
\rest{\Phi}\tau(\mathbf u) = 
c_{\tau} + \sum_{v\in A_\tau} a_{v} u_v 
+ \sum_{\rho\in B_\tau} a_{\rho} u_{\rho} 
\qquad \forall\,\, \mathbf{u}\in \sigma^{A_\tau}\times \RR^{B_\tau}_{\geq 0},
\end{equation}
for rational $c_{\tau}$ and integers $a_{v}, a_{\rho}$. Since the  $c_{\tau}$ belongs to the valuation group by assumption, we can find $\gamma_{\tau} \in K$ with $\ss \nu_{K}(\gamma_{\tau}) = c_{\tau}$. 
If we set 
\[
g_{\tau} \colonequals \gamma_{\tau}\, \prod_{v\in A}\z_{v}^{a_{v}} \prod_{\rho\in B}\w_{\rho}^{a_{\rho}} \in \Rat(X) = \Rat(\Xscr), 
\]
then $\rest{\trop(g_{\tau})}{\tau} = \rest{\Phi}\tau$, which gives the second half of the proposition.  

To get the first half, we write 
\[
g_{\tau} = \gamma_{\tau}\frac{\prod_{v\in A}\z_{v}^{\max\{a_{v}, 0\}} \prod_{\rho\in B}\w_{\rho}^{\max\{a_{\rho}, 0\}} }{\prod_{v\in A}\z_{v}^{-\min\{a_{v}, 0\}} \prod_{\rho\in B}\w_{\rho}^{-\min\{a_{\rho}, 0\}} 
}. 
\]
Let $\eta$ be a cell of $\Delta(\Xscr, \cH)$. We need to show that 
$\trop(g_{\tau})$ is integral $\QQ$-affine on $\eta$.  
By the choice of the parameters $\z_v$ and $\w_{\rho}$, if $v \notin \zeta(\eta)$ (resp. $\rho\not\in \zeta(\eta)$), then the parameter $\z_v$ (resp. $\w_\rho$) is a unit in $\OO_{\Xscr, \xi_\eta}$. Moreover, if $v \in \zeta(\eta)$ (resp. $\rho \in \zeta(\eta)$), the parameter $\z_v$ (resp. $\w_\rho$) is in $\OO_{\Xscr, \xi_\eta}$ and has order of vanishing one along the component $Z_v$ (resp. $W_\rho$). It follows that both the numerator and the  denominator are admissible expansions in the local ring $\OO_{\Xscr, \xi_\eta}$ for each $\eta \in \Delta(\Xscr, \cH)$. This implies that the tropicalization of $g_\tau$ is integral $\QQ$-affine on $\eta$. 
\end{proof}

\subsection{Proof of Theorem~\ref{thm:Q:1}}
Before giving the proof of Theorem~\ref{thm:Q:1}, we state useful lemmas. 

\begin{Lemma}
\label{lemma:generic}
Let $x_{ij}\in K$ for $1 \leq i \leq m$ and $1 \leq j \leq n$. Assume that the residue field $\kappa$ has cardinality larger than $n$. Then, there exist $\beta_i\in R^\times$ for $1 \leq i \leq m$ such that, for any $1 \leq j \leq n$, we have $|\sum_{i=1}^m \beta_i x_{ij}| = \max\{|x_{1j}|, \ldots, |x_{mj}|\}$. 
\end{Lemma}
\begin{proof}
We can assume without loss of generality that $\max\{|x_{1j}|, \ldots, |x_{mj}|\} =1$ for all $j=1, \dots, n$. Let $\bar x_{ij}$ be the reduction of $x_{ij}$ in the residue field. Consider the $m$ by $n$ matrix $M=\bigl(\bar x_{ij}\bigr)$. It will be enough to show the existence of a vector $y=(y_1, \dots, y_m)$ in $\kappa^m$ such that the coordinates of $yM$ are all nonzero. Indeed, for such $y$, we will choose $\beta_j \in R^\times$ with reduction in $\kappa$ equal to $y_j$. It will follow that the reduction of $\sum_{i=1}^m \beta_i x_{ij}$ in $\kappa$ will be nonzero, and we conclude that $|\sum_{i=1}^m \beta_i x_{ij}| = 1= \max\{|x_{1j}|, \ldots, |x_{mj}|\}$, as required.

It remains to show the existence of $y$. For each $j$, the equation $\sum_{i=1}^m y_i \bar x_{ij}=0$ defines a hyperplane in $\kappa^m$. Since the cardinality of $\kappa$ is larger than $n$, the union of these $n$ hyperplanes does not cover the full space, and therefore, there exists $y\in \kappa^m$ which does not belong to any of these hyperplanes. For this element, the coordinates of $yM$ are all nonzero.
\end{proof}

Let $\ss\Theta_{A\cup B} \colonequals 
\sigma^{A} \times {\RR^B_{\geq 0}}$ and 
$\ss\Theta_{A^\prime \cup B^\prime} \colonequals 
\sigma^{A^\prime} \times {\RR^{B^\prime}_{\geq 0}}$ be extended simplices with $A \subseteq  A^\prime$ and $B \subseteq B^{\prime}$.
We allow $A = B = \emptyset$ in which case we set 
$\ss\Theta_{A\cup B} = \ss\Theta_{\emptyset\cup \emptyset} = \emptyset$, but we assume that $A^\prime \neq \emptyset$ and so  
$\ss\Theta_{A^\prime\cup B^\prime} \neq \emptyset$.  We regard $\ss\Theta_{A\cup B}$ as the subset of $\ss\Theta_{A^\prime\cup B^\prime}$ given by
\[
\ss \Theta_{A\cup B}= 
\bigl\{\mathbf(u) \in \Theta_{A^\prime\cup B^\prime} \,\bigl|\, u_v=u_\rho =0 \textrm{ for all } v\in A^\prime \setminus A, \rho \in B^\prime\setminus B \bigr\}. 
\]

\begin{Lemma}
\label{lem:simple}
Let $F$ be a piecewise integral $\QQ$-affine function 
and $G$ an integral $\QQ$-affine function
on $\ss\Theta_{A^\prime\cup B^\prime}$ such that
\[
\rest{F}{\ss \Theta_{A\cup B}} \equiv 0, \quad \rest{G}{\ss \Theta_{A\cup B}} \equiv 0, \quad \text{and}\quad \rest{F}{\ss\Theta_{A^\prime\cup B^\prime} \setminus \ss \Theta_{A\cup B}} >0.
\]
Assume that $F$ is concave. 
For $\rho \in B^\prime \setminus B$, consider the vector $\e_\rho \in \RR^{B^\prime}$. We assume that  
\begin{itemize}
    \item[\textup{(i)}] Either, the slope of $F$ near infinity in the direction of $\e_{\rho}$ is positive, 
    \item[\textup{(ii)}] Or, if the slope of $F$ near infinity along $\e_{\rho}$ is zero, then the slope of $G$ near infinity along $\e_{\rho}$ is also zero.
\end{itemize}
Then, for any sufficiently large $N \gg 1$, we have $\rest{(N F - G)}{\ss\Theta_{A^\prime \cup B^\prime} \setminus \ss\Theta_{A \cup B}} >0$. 
\end{Lemma}

\begin{proof} 
We choose a refinement of $\Theta_{A^\prime \cup B^\prime}$ such that 
on any polyhedron of the refinement $F$ and $G$ are integral $\QQ$-affine. 
Since $F$ is strictly positive on $\Theta_{A^\prime \cup B^\prime} \setminus \Theta_{A \cup B}$, for any sufficiently large $N$, $NF-G$ is strictly positive on compact cells of the refinement.  We claim that for any half-line on an unbounded cell of the refinement, the slope of $NF-G$ is non-negative for sufficiently large $N$. Since 
$NF-G$ is concave on $\Theta_{A^\prime \cup B^\prime}$, it is enough to prove the statement for the slopes of $NF-G$ near infinity along the rays of $\Theta_{A^\prime \cup B^\prime}$. This comes from our assumptions~(i) and~(ii). To conclude, for any unbounded cell of the refinement, since $NF-G$ is affine on the cell, takes positive values on vertices away from $\Theta_{A \cup B}$, and has non-negative slopes along rays, it remains positive on the cell everywhere outside $\Theta_{A \cup B}$. 
\end{proof}

We begin the proof of Theorem~\ref{thm:Q:1}. With the notation of the previous sections,  
we set 
\begin{equation}
\label{eqn:def:h}
\phi\colonequals \sum_{\tau \in \Delta(\Xscr, \cH)}
\beta_\tau f_\tau^{N_\tau} g_\tau \in \Rat(\Xscr) = \Rat(X)
\end{equation}
as in \eqref{eqn:construction:h}, where the coefficients $\beta_\tau \in R^\times$  and the integers $N_\tau \gg 1$ will be suitably chosen later. We claim that $\trop(\phi) = \Phi$ on  $\Sk(\Xscr, \cH)$.

Fix an element $\tau \in \Delta(\Xscr, \cH)$. Then, we write 
\[
\phi =  \beta_\tau f_\tau^{N_\tau} g_\tau  + 
\sum_{
\substack{
\eta \in \Delta(\Xscr, \cH)\\
\eta\neq \tau
}}
\beta_\eta f_\eta^{N_\eta} g_\eta. 
\]
We will show that $\rest{\trop(\phi)}{\tau} = \rest{\Phi}{\tau}$. 

\smallskip

First, note that since $\beta_\tau \in R^\times$, $\rest{\trop(f_\tau)}{\tau} \equiv 0$ and $\rest{\trop(g_\tau)}{\tau} = \rest{\Phi}{\tau}$, we have 
\[
\rest{\trop(\beta_\tau f_\tau^{N_\tau} g_\tau )}{\tau}
= -\log|\beta_\tau| + N_\tau \rest{\trop(f_\tau)}{\tau} + \rest{\trop(g_\tau)}{\tau} 
= \rest{\Phi}{\tau}.
\]
We now evaluate $\rest{\trop(\beta_\eta f_\eta^{N_\eta} g_\eta)}{\tau}$ for any $\eta \neq \tau$, dividing into three cases. 

\medskip
{\bf Case 1.}\quad 
Suppose that $\eta \prec \tau$ ($\eta \neq \tau$). Then, by the property of $f_\eta$, we have $\rest{\trop(f_\eta)}{\eta} \equiv 0$ and 
$\rest{\trop(f_\eta)}{\tau\setminus \eta} > 0$. We also have $\rest{\trop(g_{\eta})}{\eta} = \rest{\Phi}{\eta}$. 

Restricting to $\eta$, we have 
\[
\rest{\trop(\beta_\eta f_\eta^{N_\eta} g_\eta)}{\eta} = \rest{\Phi}{\eta} 
= \rest{(\rest{\Phi}{\tau})}{\eta} 
= \rest{\left(\rest{\trop(\beta_\tau f_\tau^{N_\tau} g_\tau)}{\tau}\right)}{\eta}
= \rest{\trop(\beta_\tau f_\tau^{N_\tau} g_\tau)}{\eta}.
\] 
And, on $\tau \setminus \eta$, we have 
\[
\rest{\trop(\beta_\eta f_\eta^{N_\eta} g_\eta)}{\tau \setminus \eta} 
- \rest{\trop(\beta_\tau f_\tau^{N_\tau} g_\tau)}{\tau \setminus \eta} 
= \left(N_\eta \rest{\trop(f_\eta)}{\tau\setminus \eta} + \rest{\trop(g_\eta)}{\tau\setminus \eta} \right) - \rest{\Phi}{\tau \setminus \eta}. 
\]

We apply Lemma~\ref{lem:simple} to 
$\Theta_{A_\tau \cup B_\tau}$ and $ \Theta_{A_\eta\cup B_\eta}$, and $F = \rest{\trop(f_\eta)}{\tau}$ and $G = \rest{\Phi}{\tau}- \rest{\trop(g_\eta)}{\tau}$. Then, $F$ is concave, piecewise integral $\QQ$-affine and $G$ is integral $\QQ$-affine.  
We have $\rest{F}{\eta} \equiv 0$, $\rest{G}{\eta} \equiv 0$, and 
$\rest{F}{\tau \setminus \eta} > 0$. Furthermore, by the last statement of 
Proposition~\ref{prop:construction:f:sigma}, 
the slope of $F = \rest{\trop(f_\eta)}{\tau}$ near infinity along any ray of $\tau\setminus \eta$ is positive.

From Lemma~\ref{lem:simple}, we conclude that 
\[
\rest{\trop(\beta_\eta f_\eta^{N_\eta} g_\eta)}{\tau\setminus\eta} > \rest{\trop(\beta_\tau f_\tau^{N_\tau} g_\tau)}{\tau\setminus\eta} \quad\text{if}\quad N_\eta \gg 1.
\]
Here, $N_\eta$ is independent of the choice of $N_\tau$, because $\rest{\trop(f_\tau)}{\tau} \equiv 0$.

\medskip
{\bf Case 2.}\quad 
Suppose that $\eta \succ \tau$ ($\eta \neq \tau$). 
Since $\rest{\trop(f_\eta)}{\eta} \equiv 0$, we have $\rest{\trop(f_\eta)}{\tau} \equiv 0$. Also, we have
$
\rest{\trop(g_\eta)}{\tau} = \rest{\left(\rest{\Phi}{\eta}\right)}{\tau}
= \rest{\Phi}{\tau}
$. 
Thus,
\[
\rest{\trop(\beta_\eta f_\eta^{N_\eta} g_\eta)}{\tau} = \rest{\Phi}{\tau}
\quad \text{for any positive integer $N_\eta$}. 
\]

\medskip
{\bf Case 3.}\quad
In the remaining case, we have $\eta \not\preceq \tau$ and $\eta \not\succeq\tau$. We divide this into two subcases.

\smallskip

\noindent {\bf 3.1.}\quad
Suppose first that $\eta \cap \tau \neq \emptyset$. 
By the properties of $f_\eta$, we have 
\[
\rest{\trop(f_\eta)}{\tau\setminus \eta} >0\quad \text{and}\quad \rest{\trop(f_\eta)}{\tau \cap \eta} \equiv 0.
\]
Additionally, we have 
\[
\rest{\trop(g_\eta)}{\tau \cap \eta} = \rest{(\rest{\Phi}{\eta})}{\tau\cap \eta} = \rest{\Phi}{\tau\cap \eta}.
\]
Thus, we get 
\[
\rest{\trop(\beta_\eta f_\eta^{N_\eta} g_\eta)}{\tau \cap\eta} = \rest{\Phi}{\tau\cap\eta} = \rest{\trop(\beta_\tau f_\tau^{N_\tau} g_\tau)}{\tau\cap\eta}.
\]

Next, as in Case 1, apply 
Lemma~\ref{lem:simple} to 
$\ss \Theta_{A_\tau\cup B_\tau}$ and $\ss \Theta_{A_{\tau \cap \eta}\cup B_{\tau \cap \eta}}$, and two piecewise integral $\QQ$-affine functions 
$F = \rest{\trop(f_\eta)}{\tau}$ and $G = \rest{\Phi}{\tau}- \rest{\trop(g_\eta)}{\tau}$.  Again, $F$ is concave and $G$ is affine. Moreover, we have $\rest{F}{\tau \cap \eta} \equiv 0$, $\rest{G}{\tau \cap \eta} \equiv 0$, and 
$\rest{F}{\tau \setminus (\tau \cap \eta)} > 0$. 
Furthermore, by the last statement of 
Proposition~\ref{prop:construction:f:sigma}, 
the slope of $F = \rest{\trop(f_\eta)}{\tau}$ near infinity along each ray of $\tau \setminus (\tau \cap \eta)$ is positive. We conclude that 
\[
\rest{\trop(\beta_\eta f_\eta^{N_\eta} g_\eta)}{\tau\setminus\eta} > \rest{\trop(\beta_\tau f_\tau^{N_\tau} g_\tau)}{\tau\setminus\eta}
\quad \text{if $N_\eta \gg 1$}. 
\]
Note that $N_\eta$ is again independent of the choice of $N_\tau$, because $\rest{\trop(f_\tau)}{\tau} \equiv 0$.  

\smallskip

\noindent {\bf 3.2.}\quad
Suppose that $\eta \cap \tau = \emptyset$. 
In this case, we have $ \rest{\trop(f_\eta)}{\tau} = \rest{\trop(f_\eta)}{\tau\setminus \eta}>0$.

We are going to apply Lemma~\ref{lem:simple} to  $\Theta_{A_\tau\cup B_\tau}$ and $\Theta_{\emptyset \cup \emptyset}=\emptyset$, and two piecewise integral $\QQ$-affine functions $F = \rest{\trop(f_\eta)}{\tau}$ and $G = \rest{\Phi}{\tau}- \rest{\trop(g_\eta)}{\tau}$. Note that $F$ is concave and $G$ is affine.

We observe that $F > 0$. Next, we check the slope conditions for $F$ and $G$. Let $S_\tau, S_\eta \in \Delta(\Xscr, \cH)$ be the special strata corresponding to $\tau, \eta$, with generic points $\xi_\tau, \xi_\eta$,  respectively. Let $A_\tau, A_\eta \subseteq V_\f$ and $B_\tau, B_\eta \subseteq \ss V_\infty$ be the corresponding subsets. Since $\eta \cap \tau = \emptyset$, we have $A_\tau \cap A_\eta=\emptyset$.

Consider a ray $\rho\in B_\tau$. Two cases can occur depending on whether $\rho$ belongs to $B_\eta$ or not. In the latter case, $W_\rho$ does not appear in $D_\eta$. By the choice of $f_\eta$, since
$f_\eta$ vanishes on the complement $E_\eta$ of $D_\eta$, the slope of $F$ along the ray $\varrho$ given by $\e_\rho$ will be positive.

In the former case, let $\varrho^\prime$ be the ray of $\eta$ defined by $\rho$. Then, $\varrho$ and $\varrho^\prime$ are parallel. By the definition of $\Phi \in \Rat(\Delta(\Xscr, \cH))$, the slope of $\rest{\Phi}{\tau}$ along $\varrho$ 
and that of $\rest{\Phi}{\eta}$ along $\varrho^\prime$ are the same. 
Now, in the analogue equation to \eqref{eqn:pullback:Phi} for $\eta$,
\[
\rest{\Phi}{\eta}(\mathbf u) = c_\eta + \sum_{v\in A_\eta} a^\prime_v u_v + \sum_{\rho\in B_\eta} a^\prime_{\rho} u_{\rho} \qquad \forall\,\, \mathbf{u}\in \sigma^{A_\eta}\times \RR^{B_\eta}_{\geq 0},
\]
consider the coefficient $a^\prime_{\rho}$ of $u_{\rho}$. Let $\w_\rho = 0$ be the local equation for $W_\rho$ at the generic point $\xi_\eta$. By the definition of $g_\eta$, the exponent of $\w_\rho$ in $g_\eta$ is equal to $a^\prime_{\rho}$, and so the order of $g_\eta$ along the horizontal divisor $W_\rho$ is equal to $a^\prime_{\rho}$. It follows that the slope of $\rest{\trop(g_\eta)}{\tau}$ along $\varrho$ and 
that of $\rest{\trop(g_\eta)}{\eta}$ along $\varrho^\prime$ are the same integer $a^\prime_{\rho}$. Thus, for $\rho\in B_\eta \cap B_\tau$, 
the slope of $G$ along $\ss\e_{\rho}$ is zero. 

We have confirmed that $F$ and $G$ satisfy the slope conditions in Lemma~\ref{lem:simple}. We apply that lemma to conclude that 
\[
\rest{\trop(\beta_\eta f_\eta^{N_\eta} g_\eta)}{\tau\setminus\eta} > \rest{\trop(\beta_\tau f_\tau^{N_\tau} g_\tau)}{\tau\setminus\eta}
\quad \text{if $N_\eta \gg 1$}. 
\]
Note that, once again, $N_\eta$ is independent of the choice of $N_\tau$, because $\rest{\trop(f_\tau)}{\tau} \equiv 0$.

\begin{proof}[Proof of Theorem~\ref{thm:Q:1}]
In the definition of $\phi$ given in \eqref{eqn:def:h}, for each $\tau \in \Delta(\Xscr, \cH)$, we choose $N_\tau \gg 1$ sufficiently large (which depends on $\trop(g_\tau)$ and $\Phi$, but not on cells of  $\Delta(\Xscr, \cH)$ other than $\tau$).  
\smallskip

By Lemma~\ref{lemma:generic},  since $\kappa$ has large enough cardinality, we can choose $\beta_\tau \in R^\times$ suitably so that by what we have seen in Cases 1, 2, 3, we get, for any $\tau$, 
\begin{align}
\label{eqn:trop:phi:tau}
\rest{\trop(\phi)}{\tau}
&= 
\trop\left(\beta_\tau f_\tau^{N_\tau} g_\tau 
+ 
\sum_{
\substack{
\eta \in \Delta(\Xscr, \cH)\\
\eta\prec \tau
}}
\beta_\eta f_\eta^{N_\eta} g_\eta 
+ 
\sum_{
\substack{
\eta \in \Delta(\Xscr, \cH)\\
\eta \succ \tau
}}
\beta_\eta f_\eta^{N_\eta} g_\eta 
\right. 
\\
\notag
& \qquad 
\rest{
\left. 
+ 
\sum_{
\substack{
\eta \in \Delta(\Xscr, \cH)\\
\eta \not\preceq \tau,\, \eta \not\succeq \tau,\, 
\eta\cap \tau\neq \emptyset
}}
\beta_\eta f_\eta^{N_\eta} g_\eta 
+ 
\sum_{
\substack{
\eta \in \Delta(\Xscr, \cH)\\
\eta \not\preceq \tau,\, \eta \not\succeq \tau,\,
\eta\cap \tau= \emptyset
}}
\beta_\eta f_\eta^{N_\eta} g_\eta 
\right)
}{\tau} 
= \rest{\trop(\beta_\tau f_\tau^{N_\tau} g_\tau)}{\tau} 
= \rest{\Phi}{\tau}.
 \end{align} 
Indeed, recall from Case~2 that $g_\tau$ and $g_\eta$'s with $\eta \succ \tau$ are monomials whose tropicalizations give the same function $\rest{\Phi}{\tau}$ on $\tau$, and the constant term of $f_\tau$ and $f_\eta$'s  with $\eta \succ \tau$ are all units in $\widehat{\OO}_{\Xscr, \xi_\tau}$, where $\xi_\tau$ is the generic point of the special stratum corresponding to $\tau$. Thus, we can choose $\beta_\tau \in R^\times$ ($\tau \in \Delta(\Xscr, \cH)$) suitably so that 
$\rest{
\trop\left(\beta_\tau f_\tau^{N_\tau} g_\tau 
+ 
\sum_{
\substack{
\eta \in \Delta(\Xscr, \cH)\\
\eta \succ \tau
}}
\beta_\eta f_\eta^{N_\eta} g_\eta 
\right)}{\tau}
= \rest{\trop(\beta_\tau f_\tau^{N_\tau} g_\tau)}{\tau}$ holds for any $\tau$. Since $-\log |x| > -\log|y|$ implies $-\log |x+y|=-\log |y|$ for any nonzero elements $x, y$ in a non-Archimedean field, Cases~1 and ~3 then give \eqref{eqn:trop:phi:tau}. 

We infer that $\trop(\phi) = \Phi$ on $\Sk(\Xscr, \cH)$, completing the proof of Theorem~\ref{thm:Q:1}.
\end{proof}


\section{Finite generation}
\label{sec:strategy}
Let $\Delta$ be an abstract polyhedral complex with recession fan $\Upsilon$, set of vertices $V_\f$ and set of rays $\ss V_\infty$. For a finite set $C$, we denote by $\ss\e_C$ the point of $\RR^C$ with coordinates all equal to one.

For $A\subseteq V_\f$ and $B\subseteq \ss V_\infty$, let $\ss\Theta_{A\cup B} = \sigma^A\times \RR^B_{\geq 0}$ be the corresponding extended simplex in $\RR^{A\cup B} \subseteq \RR^{V_\f\cup \ss V_\infty}$. The point $p = \frac1{|A|}\ss\e_{A\cup B}= \frac1{|A|}(\ss\e_{A}, \ss \e_B)$ belongs to the extended simplex $\ss\Theta_{A\cup B}$, and we refer to it as the \emph{barycenter} of $\ss\Theta_{A\cup B}$. For each $\tau\in \Delta$ with $\zeta(\tau) = A \cup B$, we denote by $p_\tau$ the barycenter of $\tau$.

\begin{figure}[ht]
    \centering
\[
\begin{tikzpicture}[scale = 2.5]
\coordinate (v1) at (0, 1);
\coordinate (v2) at (0, 0);

\coordinate (v1end) at ($(v1) + (2,0)$);
\coordinate (v2end) at ($(v2) + (2,0)$);

\coordinate (v12) at ($1/2*(v1) + 1/2*(v2)$);
\coordinate (v1r) at ($(v1) + (1/1.4,0)$);
\coordinate (v2r) at ($(v2) + (1/1.4,0)$);
\coordinate (v12r) at ($(v12) + (1/2,0)$);
\coordinate (v12end) at ($(v12) + (2,0)$);
\coordinate (v12rbis) at ($(v12r) + (.05,0)$);

\fill[red!5] (v1) -- (v1end) -- (v2end) -- (v2) -- cycle;

\draw[thick] (v1) -- (v2);
\draw[thick] (v1) -- (v1end);
\draw[thick] (v2) -- (v2end);

\draw[red, dashed] (v1) -- (v12r) -- (v12);
\draw[red, dashed] (v2r) -- (v12r) -- (v1r);
\draw[red, dashed] (v2) -- (v12r);
\draw[red, dashed] (v12r) -- (v12end);

\draw[line width=0.1mm] (v1) circle (0.3mm);
\filldraw[aqua] (v1) circle (0.2mm);
\draw[line width=0.1mm] (v2) circle (0.3mm);
\filldraw[aqua] (v2) circle (0.2mm);
\draw[line width=0.1mm] (v12) circle (0.2mm);
\filldraw[red] (v12) circle (0.1mm);
\draw[line width=0.1mm] (v1r) circle (0.2mm);
\filldraw[red] (v1r) circle (0.1mm);
\draw[line width=0.1mm] (v2r) circle (0.2mm);
\filldraw[red] (v2r) circle (0.1mm);
\draw (v12r) circle (0.2mm);
\filldraw[red] (v12r) circle (0.1mm);

\node[below left] at (v2) {$v_2$};
\node[below left] at (v1) {$v_1$};
\node[below right] at (v2r) {$p_{2\rho}$};
\node[below right] at (v1r) {$p_{1\rho}$};
\node[below right] at (v12rbis) {$p_{12\rho}$};
\node[below left] at (v12) {$p_{12}$};
\node[right] at (v1end) {$\rho$};
\end{tikzpicture}
\]
\caption{Barycentric subdivision of an extended simplex of dimension two with two vertices and a single ray.}\label{fig:barycentric-subdivision}
\end{figure}
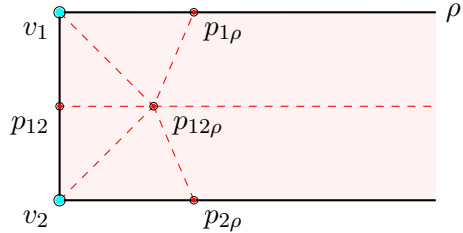

\begin{figure}[t!]
   \includegraphics{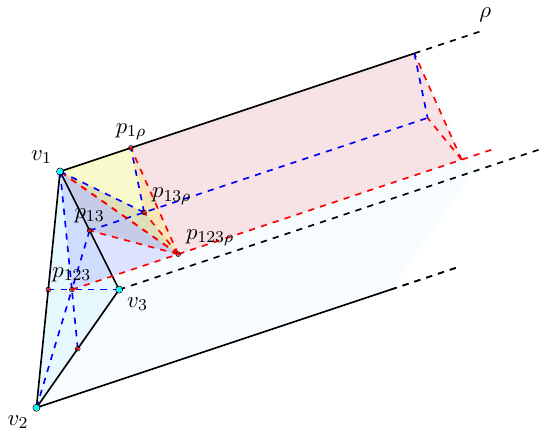}   
   \caption{Barycentric subdivision of an extended simplex with three vertices and a single ray. Only a portion of the subdivision is given, with four solids depicted in color: three tetrahedra with vertices $\{v_1, \ss p_{13}, \ss p_{123}, \ss p_{123\rho}\}$, $\{v_1, \ss p_{13}, \ss p_{13\rho}, \ss p_{123\rho}\}$, and $\{v_1, \ss p_{1\rho}, \ss p_{13\rho}, \ss p_{123\rho}\}$, and a prism with vertices $\{\ss p_{1\rho}, \ss p_{13\rho}, \ss p_{123\rho}\}$ and a ray $\rho$.}
    \label{fig:barycentric-subdivision2}
\end{figure}

 The \emph{barycentric subdivision} $\Delta'$ of $\Delta$ is defined as follows.  The vertex set $V_\f^\prime$ of $\Delta'$ is the collection of the barycenters $p_\tau$ one for each $\tau$ in $\Delta$. A chain $\tau_\bullet = (\tau_{i})_{i=1}^k$ in $\Delta$ of length $k$ is a collection of elements $\hatzero \neq \tau_1 \prec \tau_2\prec \dots \prec \tau_k$. Given a subset $B\subset \ss V_\infty$, a chain $S_\bullet$ in $B$ of length $l\in\ZZ_{\geq 0}$ is a collection of nonempty subsets $S_1\subset \dots \subset S_l$ of $B$. For $l=0$, we get the empty chain with no element. To each pair $(\tau_\bullet, S_\bullet)$ consisting of a chain $\tau_\bullet$ in $\Delta$ of some length $k \in \ZZ_{\geq 1}$ and a chain $S_\bullet$ in $\zeta(\tau_1) \cap \ss V_\infty$ of some length $l\in\ZZ_{\geq 0}$, we associate the polyhedron $\eta(\tau_\bullet, S_\bullet) \subseteq \tau_k$ defined by taking the convex-hull in $\tau_k$ of the points $p_{\tau_1},\dots, p_{\tau_k}$ and all the rays $\RR_{\geq0} \ss \e_{S_1}, \dots, \RR_{\geq 0}\ss\e_{S_l}$. Then, $\Delta^\prime$ is the collection of all the polyhedra $\eta(\tau_\bullet, S_\bullet)$, see Figures~\ref{fig:barycentric-subdivision} and~\ref{fig:barycentric-subdivision2}.

 Note that for each vertex $v \in V_\f$, we have $p_v=v$, so that we get an inclusion $V_\f \subseteq V_\f^\prime$. It is easy to check that the recession fan  $\Upsilon^\prime$ of $\Delta^\prime$ is the barycentric subdivision of the recession fan $\Upsilon$ of $\Delta$. Denote by $\ss V^{\prime}_\infty$ the set of recession rays of $\Delta^{\prime}$. We have $\ss V_\infty\subseteq \ss V_\infty^\prime$.

Denote by $\preceq$ the partial order on $\Delta^\prime$ defined by the inclusion of faces. It is described as follows. Given elements $(\tau_\bullet, S_\bullet)$ and $(\omega_\bullet, T_\bullet)$, we have $\eta(\tau_\bullet, S_\bullet)\preceq\eta(\omega_\bullet, T_\bullet)$ provided that 
\begin{itemize}
\item $\tau_\bullet$ is included in $\omega_\bullet$, that is, each element $\tau_i$ coincides with an element $\omega_{n(i)}$ in $\omega_\bullet$, and 
\item we have $S_\bullet \subseteq T_\bullet$, that is, each element $S_j$ coincides with an element $T_{m(j)}$ in $T_\bullet$.
\end{itemize}

We next embed $\Delta'$, preserving its integral structure. Let $V^\dprime_\f$ denote a collection of elements $\p_\tau$, one for each $\tau\in \Delta$. The embedding of $\Delta^\prime$ will be defined in $\RR^{V^\dprime_\f \cup \ss {V^\prime}_\infty}$, with integral structure given by the lattice $\ZZ^{V^\dprime_\f \cup \ss {V^\prime}_\infty}$.
First, we define the map
\[
\varphi \colon V^\prime_\f \to \RR^{V_\f^\dprime \cup \ss {V^\prime}_\infty}, \qquad \varphi(p_\tau) \colonequals \frac1{|\zeta(\tau)\cap V_\f|}\e_{\p_\tau}.
\]
Specifically, if $\zeta(\tau) = A \cup B$ with $A \subseteq V_\f$ and $B\subseteq \ss V_\infty$, then we set $\varphi(p_\tau) = \frac 1{|A|}\e_{\p_\tau}$. 

For each pair $(\tau_\bullet, S_\bullet)$ with $\tau_\bullet$ a chain of length $k$ in $\Delta$ and $S_\bullet$ a chain of length $l$ in $\zeta(\tau_1)\cap \ss V_\infty$, let $\sigma(\tau_\bullet, S_\bullet)$ denote the convex hull in $\RR^{V^\dprime_\f\cup \ss{V^\prime}_\infty}$ of the points $\varphi(p_{\tau_1}), \dots, \varphi(p_{\tau_k})$ and the rays $\RR_{\geq} \ss \e_{S_1}, \dots, \RR_{\geq} \ss \e_{S_l}$. Using linear interpolation, we extend $\varphi$ to a bijective map 
\[
\varphi \colon \eta(\tau_\bullet, S_\bullet) \to \sigma(\tau_\bullet, S_\bullet)
\]
for each pair $(\tau_\bullet, S_\bullet)$ as above. 

Let $\Sigma$ be the collection of the polyhedra $\sigma(\tau_\bullet, S_\bullet)$, where $\tau_\bullet$ is a chain in $\Delta$ and $S_\bullet$ is a chain in $\zeta(\tau_1)\cap \ss V_\infty$. Note that $\Sigma$ forms a subcomplex of the product $P\times \RR^{\ss {V^\prime}_\infty}$ where $P$ is the simplex in $\RR^{V^\dprime_\f}$ obtained as the convex hull of the points $\varphi(p_\tau)$ for $\tau\in \Delta$. It follows directly from this observation that $\Sigma$ is an integral $\QQ$-affine polyhedral complex in $\RR^{V^\dprime_\f\cup \ss {V^\prime}_\infty}$. 

\begin{Proposition}
\label{prop:bij:Delta:prime:to:Sigma}
Notation as above, the maps $\varphi$ glue together to define a bijective map $\varphi \colon \supp{\Delta'} \to \supp{\Sigma}$. Moreover, under this bijection, the recession fan of $\Delta'$ is identified with the recession fan of $\Sigma$.

The maps $\varphi$ and its inverse $\varphi^{-1}\colon \supp{\Sigma}\to \supp{\Delta'}$ preserve the integral affine structures. In particular, we have $\ss\Rat(\supp{\Delta'}) = \ss\Rat(\supp{\Sigma})$.
\end{Proposition}
\begin{proof} The first assertion is straightforward. To show that $\varphi$ and its inverse preserve the integral structure, it suffices to prove that each map 
\[
\varphi \colon \eta(\tau_\bullet, S_\bullet) \to \sigma(\tau_\bullet, S_\bullet)
\]
is induced from an isomorphism between the lattice of $\eta(\tau_\bullet, S_\bullet)$ and the lattice of $\sigma(\tau_\bullet, S_\bullet)$.
 Let $N$ and $N^\prime$ denote the former and the latter lattices, respectively. 
 
Write $\tau_\bullet =\tau_1\prec \dots \prec \tau_k$, with $\zeta(\tau_i) = A_i \cup B_i$, where $A_i \subseteq V_\f$ and $B_i \subseteq \ss V_\infty$. Write $S_\bullet = S_1 \subset \dots \subset S_l\subseteq B_1=\zeta(\tau_1)\cap \ss V_\infty$. 

Let $m_{i,j}$ be the least common multiple of $|A_i|$ and $|A_{j}|$ for each pair $i,j=1, \dots, k$.  The lattice $N^\prime$ is generated by the vectors 
\[
m_{i,j}\left(\frac{1}{|A_i|}\e_{\p_{\ss \tau_i}} - \frac{1}{|A_j|}\e_{\p_{\ss \tau_j}}\right), \quad i,j=1, \dots, k, \quad \text{and}\quad \ss\e_{S_j}, \quad j=1,\dots, l.
\]
Similarly, the lattice $N$ is generated by the vectors 
\[
m_{i,j}\left(\frac{1}{|A_i|}\e_{A_i} - \frac{1}{|A_j|}\e_{A_j}\right), \quad i,j=1, \dots, k, \quad \text{and} \quad \ss\e_{S_j}, \quad j=1,\dots, l.
\]
The map $\varphi$ is the restriction of a linear map $\psi$ from the vector space spanned by  $\e_{A_i}$ and $\ss\e_{S_j}$ to the vector space spanned by $\e_{\ss\tau_i}$ and $\ss\e_{S_j}$. This map is defined by $\psi(\e_{A_i}) = \e_{\tau_i}$ and $\psi(\ss\e_{S_j})=\ss\e_{S_j}$. It is straightforward to verify that $\psi$ and its inverse are both integral maps. This establishes the preservation of the integral structure and proves the second statement in the proposition. The final assertion follows immediately as a consequence of the first two statements.
\end{proof}

Note that the embedding $\varphi$ of $\Delta^\prime$ induces an embedding of the recession fan of $\Delta^\prime$ as a subfan of the fan corresponding to the positive orthant $\RR^{\ss {V^\prime}_\infty}_{\geq 0}$. Denote by $\Upsilon$ the recession fan of $\Sigma$.
\begin{Proposition}\label{prop:completion}
There exists an integral $\QQ$-affine polyhedral complex structure $\widehat\Sigma$ on $\RR^{V_\f^\dprime\cup \ss {V^\prime}_\infty}$ that contains $\Sigma$ as a polyhedral subcomplex, and has a recession fan $\widehat \Upsilon$ which is a completion of $\Upsilon$ to a complete fan in $\RR^{V_\f^\dprime\cup \ss {V^\prime}_\infty}$. 
Additionally, we can choose $\widehat\Sigma$ to be convex with a convex recession fan $\widehat \Upsilon$.
\end{Proposition}
\begin{proof}
As observed, $\Sigma$ is a polyhedral subcomplex of the product $P\times \RR^{\ss {V^\prime}_\infty}_{\geq 0}$, where $P$ is the simplex  with vertices $\varphi(p_\tau)$ for $\tau \in \Delta$. We extend $P$ to a convex, integral $\QQ$-affine polyhedral complex structure on the entire space $\RR^{V_\f^\dprime}$, ensuring that it has a convex recession fan. 

Similarly, the positive orthant $\RR^{\ss {V^\prime}_\infty}_{\geq 0}$ is a polyhedral subcomplex of the polyhedral complex formed by all the orthants, which is also convex. The product of these polyhedral complex structures on $\RR^{V_\f^\dprime}$ and $\RR^{\ss {V^\prime}_\infty}$ yields a polyhedral complex structure $\widehat\Sigma$ with the desired properties stated in the proposition. 
\end{proof}

We fix a polyhedral tubular neighborhood $U$ of $\supp{\Sigma}$ in $\RR^{V_\f^\dprime \cup \ss {V^\prime}_\infty}$, which deformation retracts to $\supp{\Sigma}$ via a facewise affine map.  

Next, we select a family of integral $\QQ$-affine polyhedral complexes $(\widehat\Sigma_\alpha)_\alpha$ with support $\RR^{V^\dprime_\f\cup \ss {V^\prime}_\infty}$, such that:
\begin{itemize}
    \item Each $\widehat \Sigma_\alpha$ is a simplicial refinement of $\widehat \Sigma$;
    \item For each $\alpha\in \NN$, there exists a strictly convex function $\widehat H_\alpha \in \Rat(\widehat\Sigma_\alpha)$; and
    \item For any refinement $\widehat\Sigma^\prime$ of $\widehat\Sigma$,  there exists an $\alpha$ such that $\widehat\Sigma_\alpha$ is a refinement of $\widehat\Sigma^\prime$.
\end{itemize}
Furthermore, by restricting to a suitable subfamily if necessary, we can assume without loss of generality that every face $\tau$ of $\widehat\Sigma_\alpha$ intersecting $\supp{\Sigma}$ lies entirely within $U$. 

For each $\alpha$, let $\Sigma_\alpha$ denote  the induced subdivision of $\supp{\Sigma}$ from $\widehat\Sigma_\alpha$. Thus, we have $\supp{\Sigma_\alpha} = \supp{\Sigma}$ and 
\[
\bigcup_\alpha \Rat(\Sigma_\alpha) = \Rat(\supp{\Sigma}).
\]
We define $\Sigma^{\prime}_\alpha$ as the union of all faces of $\widehat\Sigma_\alpha$ that have nonempty intersection with $\supp{\Sigma}$. It follows from the assumption made above that $\supp{\Sigma^{\prime}_\alpha} \subseteq U$.

\begin{Proposition} \label{prop:difference-convex}
Each facewise affine function $F$ on $\widehat \Sigma_\alpha$ can be expressed as the difference of two strictly convex facewise affine functions $H_1$ and $H_2$ on $\widehat \Sigma_\alpha$. Moreover, if $\rest{F}{\Sigma^\prime_\alpha}\in \Rat(\Sigma^\prime_\alpha)$, we can choose $H_1$ and $H_2$ such that both  $\rest{H_1}{\Sigma^\prime_\alpha}$ and $\rest{H_2}{\Sigma^\prime_\alpha}$ belong to $\ss \Rat(\Sigma^\prime_\alpha)$.
\end{Proposition}
\begin{proof}
For sufficiently large $\lambda \in \NN$, the functions $H_1 = \lambda \widehat H_{\alpha} + F$ and $ H_2=\lambda \widehat H_{\alpha}$ are both strictly convex. This proves the first part of the proposition.

To prove the second part, we observe that since $\widehat H_\alpha \in \Rat(\widehat \Sigma_\alpha)$, the restriction $\rest{\widehat H_\alpha}{\Sigma^\prime_\alpha}$ belongs to $\Rat(\Sigma^\prime_\alpha)$. Therefore, we conclude that both $H_1$ and $H_2$ are elements of $\Rat(\Sigma^\prime_\alpha)$, as required.
\end{proof}

Let $H \in \ss\Rat(\widehat \Sigma_\alpha, \RR)$ be any strictly convex facewise affine function on $\widehat\Sigma_\alpha$ such that the restriction $\rest{H}{\Sigma^\prime_\alpha}$ belongs to $\Rat(\Sigma^\prime_\alpha)$. For each $\sigma \in \widehat\Sigma_\alpha$, we choose an  affine function $L_\sigma$ such that $\rest{(H-L_\sigma)}{\sigma}=0$, and $\rest{(H-L_\sigma)}{\eta \setminus \sigma}>0$ for all $\eta \supseteq \sigma$. Moreover, since $\rest{H}{\Sigma^\prime_\alpha}\in\Rat(\Sigma^\prime_\alpha)$, we can assume that $L_\sigma$ is integral $\QQ$-affine for all $\sigma \in \Sigma_\alpha$. 

\begin{Proposition} \label{prop:maximum-affine-functions} Notation as above, we have
\[
H = \max_{\sigma\in \widehat\Sigma_\alpha} L_\sigma \qquad \text{and} \qquad \rest{H}{\Sigma_\alpha} = \max_{\sigma\in \Sigma_\alpha} \rest{L_\sigma}{\Sigma_\alpha}.
\]
\end{Proposition}
\begin{proof}
The inequality $H \leq \max_{\sigma \in \widehat\Sigma_\alpha} L_\sigma$ is straightforward. Indeed, 
let $x \in \RR^{V^\dprime_\f\cup \ss{V^\prime}_\infty}$, and suppose $x \in \tau$ for $\tau \in \widehat\Sigma_\alpha$. Then, 
\[
H(x) = L_{\tau}(x) \leq \max_{\sigma \in \widehat\Sigma_\alpha} L_\sigma(x).
\]
Similarly, we have $\rest{H}{\Sigma_\alpha} \leq \max_{\sigma \in \Sigma_\alpha} \rest{L_\sigma}{\Sigma_\alpha}$.

For the reverse inequality, we use the lemma below which asserts that a facewise affine function $H$ on $\widehat\Sigma_\alpha$ that is convex according to Definition~\ref{def:convex:function} is also convex as a real-valued function on $\RR^{V_\f^\dprime\cup \ss {V^\prime}_\infty}$. We claim that for each $\sigma \in \widehat \Sigma_\alpha$, we have $H \geq L_\sigma$ globally on $\RR^{V^\dprime_\f\cup \ss {V^\prime}_\infty}$.

Fix a point $p$ in the relative interior of $\sigma$. Let $x$ be any point of $\RR^{V^\dprime_\f\cup \ss{V^\prime}_\infty}$, and consider the segment $[p,x]$. The restriction $\rest{H}{[p,x]}$ is convex, and the restriction of $L_\sigma$ to $[p,x]$ is affine. There exists a small segment $I \subseteq [p,x]$ containing $p$ such that $\rest{H}{I} \geq \rest{L_\sigma}{I}$, with $H(p) =L_\sigma(p)$. The function $\rest{(H-L_\sigma)}{[p,x]}$ is convex, it takes value zero at $p$, and is non-negative on $I$. We conclude that $\rest{(H-L_\sigma)}{[p,x]}$ is non-negative everywhere, i.e., $\rest{H}{[p,x]} \geq \rest{L_\sigma}{[p,x]}$. This gives $H(x) \geq L_{\sigma}(x)$. 

Thus, we obtain both inequalities $H \geq \max_{\sigma \in \widehat\Sigma_\alpha} L_\sigma$ and $\rest{H}{\Sigma_\alpha} \geq \max_{\sigma \in \Sigma_\alpha} \rest{L_\sigma}{\Sigma_\alpha}$.
\end{proof}
\begin{Lemma}
 Let $\Sigma$ be a polyhedral complex with support $\RR^n$. Let $H \colon \RR^n \to \RR$ be a facewise affine function on $\Sigma$. The following properties are equivalent:
 
 \begin{enumerate}[label=\textup{(\arabic*)}]
    \item \label{equivalence_convexity1} $H$ is convex according to Definition~\ref{def:convex:function}.
    \item \label{equivalence_convexity2} $H$ is convex as a real-valued function on $\RR^n$.
 \end{enumerate}   
\end{Lemma}

\begin{proof} The result is well-known, and we provide only a sketch of the proof.  

To show \ref{equivalence_convexity1} implies \ref{equivalence_convexity2}, take a segment $[p,q]$ in $\RR^n$. It will be enough to show that $F=\rest{H}{[p,q]}$ is convex. Since $F$ is a piecewise affine function, Definition~\ref{def:convex:function} ensures that at points of $[p,q]$ where there is a change of slopes, the sum of outgoing slopes is positive, which proves the convexity. 

    To prove the other direction, consider the upper-graph $\Gamma_H^+$ of $H$, defined as the set of all points $(x,t) \in \RR^n\times \RR$ with $t\geq H(x)$. Since $H$ is convex, $\Gamma_H^+$ is a convex subset of $\RR^{n+1}$. 
    
    Let $p$ be a point in the relative interior of a face $\sigma$ of $\Sigma$. There exists an affine hyperplane $A$ passing through $(p,H(p))$ such that the upper-graph of $H$ is entirely above $A$. Since $\rest{H}{\sigma}$ is affine, we can assume that $A$ contains the affine space generated by $\sigma$. Thus, $A$ is of the form $(x, L_\sigma(x))$ for an affine function $L_\sigma$ on $\RR^n$. We have $\rest{H}{\sigma} = \rest{L_\sigma}{\sigma}$, and the inequalities $(H-L_\sigma)|_{\eta\setminus \sigma} \geq 0$ hold for all $\eta \supset \sigma$, as required. 
\end{proof}    

\subsection{Proof of Theorem~\ref{thm:Q:2}}  Taking the barycentric subdivision $\Delta^\prime$ of $\Delta$ as above, we embed $\Delta^\prime$ as an integral $\QQ$-affine polyhedral complex $\Sigma$ in some linear space $\RR^n=\RR^{V^\dprime_\f\cup \ss{V^\prime}_\infty}$, with the notation introduced above. We also keep the notation for $\widehat \Sigma$, $(\widehat\Sigma_\alpha)_\alpha$, $\Sigma^\prime_\alpha$. Let $\pi \colon U \to \supp{\Sigma}$ be the fixed retraction map induced from the polyhedral tubular neighborhood $U$ of $\supp{\Sigma}$ to $\supp{\Sigma}$.
 
Note that we have $\Rat(\supp{\Delta})= \Rat(\supp{\Sigma})=\Rat(\supp{\Sigma_\alpha})$ for all $\alpha$. We prove $\Rat(\supp{\Sigma})$ is generated over $\TT\QQ$ by the restriction of coordinate functions $\x_1, \dots, \x_n$ of $\RR^n$ to $\supp{\Sigma}$.

 Let $F\in \Rat(\supp{\Sigma})$ be a piecewise integral $\QQ$-affine function. There exists some $\alpha$ such that $F$ belongs to $\Rat(\Sigma_\alpha)$. Using the retraction map $\rest{\pi}{\Sigma^\prime_\alpha} \colon \supp{\Sigma^\prime_\alpha} \to \supp{\Sigma}$, we extend $F$ to an integral facewise $\QQ$-affine function $G$ on $\Sigma^\prime_\alpha$, i.e., we set $G=F\circ \rest{\pi}{\Sigma^\prime_\alpha} \in \ss\Rat(\Sigma^\prime_\alpha)$.

Let $\widehat \Upsilon_\alpha$ denote the recession fan of $\widehat \Sigma_\alpha$. 

Next, we extend $G$ to a facewise (not necessarily integral $\QQ$-) affine function $\widehat G$ on $\widehat\Sigma_\alpha$. The extension process works as follows.
 
 First, for each ray $\varrho$ in $\widehat \Upsilon_\alpha$, if $\varrho$ is parallel  to a ray $\rho$ in $\Upsilon_\alpha$, we denote by $s_\varrho$ the slope of $F$ along any ray  of $\Sigma$ parallel to $\rho$. This is well-defined because $F \in \Rat(\Sigma_\alpha)$. 
 
 Then, for rays $\varrho$ of $\widehat \Upsilon_\alpha$ which are not parallel to any ray in $\Upsilon_\alpha$, we set $s_\varrho =0$. 
 
 We define $\widehat G(v)=0$ for every vertex $v$ of $\widehat\Sigma_\alpha$ not contained in $\Sigma^\prime_\alpha$. 
 
 Once the values of $\widehat G$ at the vertices of each polyhedron $\tau$ in $\widehat\Sigma_\alpha$ and the slopes along the rays of $\tau$ are fixed, since $\tau$ is simplicial, we can extend $\widehat G$ by linear interpolation to the entire polyhedron $\tau$. This way, we obtain an element $\widehat G\in \Rat(\widehat \Sigma_\alpha, \RR)$. Note that the restriction $\rest{\widehat G}{\Sigma^\prime_\alpha}$ is integral $\QQ$-affine.

 By Proposition~\ref{prop:difference-convex}, we can write $\widehat G = H_1-H_2$ for strictly convex functions $H_1$ and $H_2$ on $\widehat\Sigma_\alpha$ such that the restrictions $\rest{H_1}{\Sigma^\prime_\alpha}$ and $\rest{H_2}{\Sigma^\prime_\alpha}$ both belong to $\Rat(\Sigma^\prime_\alpha)$. By Proposition~\ref{prop:maximum-affine-functions}, we can express $\rest{H_1}{\Sigma_\alpha}$ and $\rest{H_2}{\Sigma_\alpha}$ as
 \[
 \rest{H_1}{\Sigma_\alpha} =\max_{\sigma\in \Sigma_\alpha} \rest{L_{1,\sigma}}{\Sigma_\alpha} \qquad \text{and} \qquad \rest{H_2}{\Sigma_\alpha} =\max_{\sigma\in \Sigma_\alpha} \rest{L_{2,\sigma}}{\Sigma_\alpha}
 \]
 where $L_{1,\sigma}, L_{2,\sigma}$ are integral $\QQ$-affine functions on $\RR^n$ with 
 \[
  \rest{H_1}{\sigma}=\rest{L_{1,\sigma}}{\sigma} \quad \text{and} \quad \rest{H_2}{\sigma}=\rest{L_{2,\sigma}}{\sigma}
 \] 
for all $\sigma \in \Sigma_\alpha$. 

From the above, we conclude that $F$ belongs to the semifield generated by $\TT\QQ$ and the restriction of the coordinate functions $\x_1, \dots, \x_n$ of $\RR^n$ to $\supp{\Sigma}$. Thus, the theorem is proven. \qed

\subsection{Proof of Theorem~\ref{thm:Q:2:stronger}}
To prove the first statement, using the work by Ewald and Ishida~\cite{EI06}, we extend $\Sigma$ to an integral $\QQ$-affine polyhedral complex $\widehat\Sigma$ with full support $\RR^n$ which verifies the properties stated in Proposition~\ref{prop:completion}. The rest of the arguments leading to the proof of Theorem~\ref{thm:Q:2} remains unchanged and gives the first part of Theorem~\ref{thm:Q:2:stronger}.

We prove the second part. For $W \subseteq \RR^n$ a subset of $\RR^n$, let
\begin{align*}
E(W) = \left\{(F, G) \in \Rat(\RR^n)^2  \mid 
\text{$F(w) = G(w)$ for all $w \in W$}\right\}. 
\end{align*}
Then, $E(W)$ is a congruence on $\Rat(\RR^n)$. It is proved in \cite{Song} that the following are equivalent:
\begin{enumerate}
\item[\textup{(i)}]
$E(W)$ is a finitely generated congruence.
\item[\textup{(ii)}]
$W$ is the support of an integral $\QQ$-affine polyhedral complex. 
\end{enumerate}
(The result is stated over $\TT$ in \emph{loc.cit.}, but the same proof works over $\TT\QQ$.) 

Consider the map
$\pi\colon \Rat(\RR^n) \to \ss\Rat(\supp{\Sigma})$ from the introduction. The kernel of $\pi$ is by definition equal to $E(\supp{\Sigma})$. Since $\supp{\Sigma}$ verifies (ii), $E(\supp{\Sigma})={\rm Ker}(\pi)$ is a finitely generated congruence. \qed

\subsection{Proof of Theorem~\ref{thm:GRW}}
Let $X$ be a smooth projective variety over $K$ 
with a strictly semistable pair $(\Xscr, \cH)$. As before, 
we denote by $\Delta(\Xscr, \cH)$ the dual complex of $(\Xscr, \cH)$ with the set of vertices~$V_\f$ and the set of rays $\ss V_\infty$, and we identify $\supp{\Delta(\Xscr, \cH)}$ with $\Sk(\Xscr, \cH)$. 

Let 
\[
\varphi\colon \supp{\Delta(\Xscr, \cH)^\prime} \to \supp{\Sigma} \subseteq  \RR^{V_\f^{\dprime} \cup V_\infty}
\]
be the map in Proposition~\ref{prop:bij:Delta:prime:to:Sigma}, 
so that $\varphi$ is bijective and that 
both $\varphi$ and $\varphi^{-1}$ preserve the integral $\QQ$-affine structures. 
Here, $\Delta(\Xscr, \cH)^\prime$ is a subdivision of $\Delta(\Xscr, \cH)$ with the same recession fan, 
and we have 
$\supp{\Delta(\Xscr, \cH)^\prime} = 
\supp{\Delta(\Xscr, \cH)} = \Sk(\Xscr, \cH)$.

Let $\x_1, \ldots, \x_n$ be the coordinate functions on $\RR^{V_\f^{\dprime} \cup V_\infty}$, where $n = \dim \RR^{V_\f^{\dprime} \cup V_\infty}$. We define $F_i \colonequals \varphi^* \x_i \in \ss\Rat(\supp{\Delta(\Xscr, \cH)^\prime})$. Then we have $\varphi = (F_1, \ldots, F_n)$. 
Since $\ss\Rat(\supp{\Sigma})$ is generated by 
$\x_1, \ldots, \x_n$ over $\TT\QQ$ by Theorem~\ref{thm:Q:2}, 
$\ss\Rat(\supp{\Delta(\Xscr, \cH)^\prime})$ is 
generated by $F_1, \ldots, F_n$ over $\TT\QQ$. 

By Theorem~\ref{thm:Q:1}, we take $f_i \in \Rat(X_{\bar{K}})$ such that $\trop(f_i) = F_i$. 
Then setting $U = X_{\bar{K}} \setminus \bigcup_{i=1}^n \supp{\zero(f_i)}$, the map 
\[
\psi\colon 
U^{\an} \to \RR, \quad x \mapsto (-\log|f_1(x)|, \ldots, 
-\log|f_n(x)|) = (\nu_x(f_1), \ldots, \nu_x(f_n))
\]
gives a faithful tropicalization of the skeleton $\Sk(\Xscr, \cH)$. This completes the proof of Theorem~\ref{thm:GRW}. \qed

\bibliographystyle{alpha}
\bibliography{bibliography}

\end{document}